\documentclass[aap,nameyear,MSNbibl,dvips]{arximspdf}
\usepackage{graphicx}

\doi{10.1214/10-AAP724}
\volume{21}
\issue{3}
\pubyear{2011}
\firstpage{801}
\lastpage{842}

\makeatletter

\newtheorem{theorem}{Theorem}
\newtheorem{lem}[theorem]{Lemma}
\newtheorem{cor}[theorem]{Corollary}

\newcommand{\ball}{\operatorname{ball}}
\newcommand{\dd}{{d}}
\newcommand{\dist}{\operatorname{dist}}
\newcommand{\origin}{\mathbf{o}}

\newcommand{\iint}{\int\!\!\int}

\makeatother

\begin{document}
\begin{frontmatter}

\title{Geodesics and flows in a Poissonian city}
\runtitle{Poissonian city}

\begin{aug}
\author[A]{\fnms{Wilfrid S.} \snm{Kendall}\corref{}\ead[label=e1]{w.s.kendall@warwick.ac.uk}}
\runauthor{W. S. Kendall}
\affiliation{University of Warwick}
\address[A]{Department of Statistics\\
University of Warwick\\
Coventry CV4 7AL\\
United Kingdom\\
\printead{e1}} 
\end{aug}

\received{\smonth{10} \syear{2009}}
\revised{\smonth{6} \syear{2010}}

%
\begin{abstract}
The stationary isotropic Poisson line process was used to derive upper
bounds on mean excess network geodesic length in Aldous and Kendall
[\textit{Adv. in Appl. Probab.} \textbf{40} (2008) 1--21]. The current paper
presents a study of the geometry and fluctuations of near-geodesics in
the network generated by the line process. The notion of a ``Poissonian
city'' is introduced, in which connections between pairs of nodes are
made using simple ``no-overshoot'' paths based on the Poisson line
process. Asymptotics for geometric features and random variation in
length are computed for such near-geodesic paths; it is shown that they
traverse the network with an order of efficiency comparable to that of
true network geodesics. Mean characteristics and limiting behavior at
the center are computed for a natural network flow. Comparisons are
drawn with similar network flows in a city based on a comparable
rectilinear grid. A~concluding section discusses several open problems.
\end{abstract}

%

\pdfkeywords{60D05, 90B15, Dufresne integral, frustrated optimization,
geometric spanner network, growth process, improper anisotropic Poisson line
process, Lamperti transformation, Laplace exponent, Levy process,
logarithmic excess, Manhattan city network,
Mills ratio, mark distribution, martingale central limit theorem,
network geodesic, Palm distribution, perpetuity,
Poisson line process, Poissonian city network,
Slivynak theorem, spanner, spatial network,
subordinator, traffic flow, uniform integrability}
\vspace*{-3pt}
\end{frontmatter}

\section{Introduction}\label{sec:introduction}


The ``Poissonian city'' is a random network of connections based on a
Poisson line process.
\citet{AldousKendall-2007} used such a network to address a
problem in \textit{frustrated optimization}: construct planar networks
connecting a large number of nodes such that:
\begin{enumerate}
\item
the total connection length is not much
larger than the minimum possible connection length, but also such that
\item
the average connection distance between
two randomly chosen nodes is not greatly in excess of the Euclidean distance.
\end{enumerate}
It transpires that networks satisfying criterion 1 may be augmented by
sparse Poisson line processes so
as to satisfy criterion 2 as well. More\vadjust{\goodbreak}
precisely, suppose that
$n$ nodes are distributed in an arbitrary fashion (deterministically or
randomly) over a square of total area $n$.
Recall that the minimum total length for a connecting network is
achieved by a Steiner minimum tree [\citet{PromelSteger-2002}
survey Steiner trees in general; for probabilistic aspects, see
\citet{Steele-1997}, \citet{Yukich-1998}].
It is shown by
\citeauthor{AldousKendall-2007} [(\citeyear{AldousKendall-2007}),
Theorem 1(b)] that
augmentation by a sparse Poisson line process can convert a Steiner
minimum tree into a network whose total connection length is only
slightly increased but which now
delivers a mean connection distance that is no more than $O(\log n)$ in
excess of the Euclidean distance.
(Here, the mean involves averaging the choice of nodes rather than the
randomness of the line process.)
Under a suitable weak uniformity condition on the empirical spatial
distribution of the nodes, \citeauthor{AldousKendall-2007} [(\citeyear
{AldousKendall-2007}), Theorem 2]
also establish a lower bound on the mean excess: it must be of order at
least $\Omega(\sqrt{\log n})$.

The primary motivation of the previously mentioned work was to gain a
better understanding of the behavior of network statistics (such as the
mean excess network length) for entirely general networks.
However, the appearance of Poisson line processes in the upper bound
result motivates
the following, more detailed, study of the ``Poissonian city''
generated by a unit intensity stationary isotropic Poisson line process.
What can be said about the ``near-geodesics'' used to establish the
upper bound? How close are they to true network geodesics? How does
random fluctuation affect excess length? And what about traffic flow on
such a network?
These questions are addressed below; their answers require the use of
L\'evy subordinators, self-similar Markov processes (something of a
novelty in stochastic geometry) and a curious improper anisotropic
Poisson line process.

Previous relevant work includes: the note by \citet
{Davidson-1974c}, who gives a qualitative argument showing that the
Poisson line process provides good connections;
\citeauthor{RenyiSulanke-1968} [(\citeyear{RenyiSulanke-1968}), Satz
5], who derive a result similar to the mean-excess
result, but concerning numbers of edges rather than length, and based
on a fixed number of random lines; and recent higher-dimensional
generalizations of the R\'enyi--Sulanke work by \citet
{BoroczkySchneider-2008},
Theorem~1.3. We also mention work by \citet
{VossGloaguenSchmidt-2009} on limit distributions of shortest paths
from subsidiary to major nodes in hierarchical networks based on random
tessellations. Finally, we note the interesting work of \citet
{BaccelliTchoumatchenkoZuyev-2000} related to the concept of
\textit{spanners} from graph theory [a geometric spanner is a planar graph
connecting a set of nodes for which the graph distance between any two
points is less than some fixed multiple of Euclidean distance; see,
e.g., the exposition \citet{NarasimhanSmid-2007}]. The networks
constructed in \citet{AldousKendall-2007} are averaged rather
than uniform versions of geometric spanner networks, for which the
fixed multiple of Euclidean distance is replaced by a logarithmic
additive excess and a specific constraint is imposed on the total
network length (rather than, say, small vertex degree or total number
of edges).

In the remainder of this introductory section we introduce basic
notation and concepts, and enumerate the questions to be addressed
concerning the behavior of near-geodesics and traffic flow in the
Poissonian city.

\subsection{Notation and basic concepts}\label{sec:caricature}
We begin by presenting a brief summary of stationary isotropic Poisson
line processes so as to fix notation and collect some facts about line
processes which will be used below. Further information can be found
in, for example, \citet{StoyanKendallMecke-1995}.
The ensemble of undirected lines $\ell$ in the plane may be viewed as
a M\"obius strip of infinite width or as a once-punctured projective
plane (since such lines can be produced as intersections of planes
through the origin in $3$-space with the plane $z=1$, in which case the
plane through the origin and parallel to $z=1$ does not produce an
intersection). It is often convenient to parametrize this ensemble of
lines $\ell$ by representing lines using points $(r,\theta)$, where
$r$ is the perpendicular signed distance from the line $\ell$ to a
reference point, and $\theta\in[0,\pi)$ is the angle that $\ell$
makes with a reference line running through the reference point.
A \textit{unit intensity stationary isotropic Poisson line process}
(``Poisson line process'' for short) is determined as a Poisson point
process on the representing projective plane using the intensity
measure $\frac{1}{2}\,{\dd} r\,\dd\theta$. The factor $\frac
{1}{2}$ ensures that the mean number of Poisson lines hitting a line
segment is equal to the length of the segment.

Slivynak's theorem on the Palm distribution of a Poisson process
applies here: if we condition on a specific line $\ell$ belonging to
the Poisson line process, then the residual line process is still unit
intensity stationary isotropic Poisson.

An alternative parametrization $(x,\theta)$, sometimes of use, employs
the line angle $\theta$ as above, with $x$ being the signed distance
along the reference line from the reference point to the intersection
of $\ell$ with the reference line. This representation breaks down
when $\theta=0$ (not a major issue in the case of a Poisson line
process, for which the set of lines at $\theta=0$ has zero
probability). In these coordinates the unit intensity measure is
$\frac{1}{2}\sin\theta\,\dd x\,\dd\theta$; the
sine-weighting corresponds to a length-biasing phenomenon when sampling
Poisson lines according to their intersections with a test line. In
particular, if two lines are conditioned to pass through a given point,
then they form an exchangeable pair: one having uniform direction, and
the angle $\alpha\in[0,\pi)$ between them having density $\frac
{1}{2}\sin\alpha$, independent of the direction of the first.\looseness=-1

Viewed as a random measure, the Poisson line process generates a
measure via the mean total length of lines intersected with a given
set. Testing against a unit disc, we can compute the resulting length
intensity as $\frac{\pi}{2}$. It follows from the above and
from Slivynak's theorem that the point process of intersections of
lines from the Poisson line process has intensity $\frac{\pi}{2}$
[\citet{Miles-1964a}, Theorem 2].\vadjust{\goodbreak}

\begin{figure}

\includegraphics{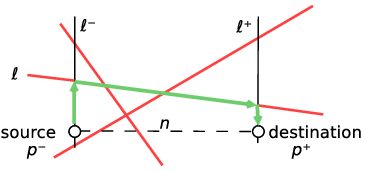}

\caption{A caricature of the procedure
of finding a route using a Poisson line process; consider only those
routes involving the two vertical lines $\ell^\pm$ and \textit{one} of
the other lines. Here a possible route is indicated by
arrows.}\label{fig:caricature}
\end{figure}

The following caricature supplies a good intuition as to where the
logarithmic excess might be located on a typical route on a network
based on a Poisson line process. Consider a network formed between just
two nodes, $p^-$ (the source) and $p^+$ (the destination), with lines
provided by a unit rate stationary isotropic Poisson line process $\Pi
$. Let the two nodes be separated by distance $n$. We condition the
line process to contain two lines $\ell^\pm$ running through source
and destination nodes, which are both perpendicular to the segment
connecting $p^-$ to $p^+$ (Figure~\ref{fig:caricature}). We consider
only those routes which involve the conditioned lines $\ell^-$
(resp., $\ell^+$) running through $p^-$ (resp., $p^+$), along with just
one other line of the Poisson line process.

%

\begin{figure}[b]

\includegraphics{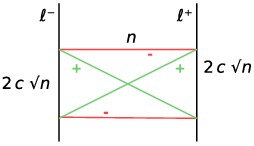}

\caption{Illustration of a classic
stochastic geometry construction for calculating the invariant measure
of the set of lines hitting both of the two vertical line segments
$\ell^\pm$ marked out by the two horizontal
lines.}\label{fig:quad}
\end{figure}

Consider the set of lines $\ell$ which intersect both $\ell^-$ at
distance at most $c\sqrt{n}$ from $p^-$ and $\ell^+$ at distance at
most $c\sqrt{n}$ from $p^+$.
Classic stochastic geometry arguments [based on inclusion--exclusion
and a special case of the ``Buffon--Sylvester problem''; see, e.g.,
\citet{Ambartzumian-1990}] then show that the invariant measure
of this line set is given by half the difference between the summed
length of the two diagonal lines minus the summed length of the two
vertical lines in Figure~\ref{fig:quad}.

%
%

Hence, the probability of \textit{no} unconditioned Poisson lines
falling in this set is
\[
\exp\bigl(
-\tfrac{1}{2}\bigl(2\sqrt{4 c^2 n + n^2} - 2n\bigr)
\bigr) \geq\exp(-2c^2) ,
\]
and, as a consequence, the resulting mean excess is bounded below by
\[
\sqrt{n}\int_0^\infty e^{-2c^2}\,\dd c
=
\frac{1}{2}\sqrt{\frac{\pi n}{2}} ,
\]
attributable to the parts of the route which lie on the conditioned
lines $\ell^-$, $\ell^+$.
(Excess along the unconditioned line $\ell$ itself is bounded above by
$\sqrt{4c^2 n +n^2}-n\leq2c^2$ and is hence negligible in the case of
large $n$.)

This rather trivial example makes it clear that being permitted to use
more than one line (in addition to the two conditioned lines) will
reduce the excess principally by rounding off the corners at the start
and finish of the journey.
Thus, it is clear [as is indeed apparent from
details of the computations in \citet{AldousKendall-2007},
Theorem 3] that the logarithmic excess in the full construction
is a cost which arises entirely from the business of getting on and off
an efficient route between source and destination.

\subsection{Making connections}\label{sec:connections}
\begin{figure}[b]

\includegraphics{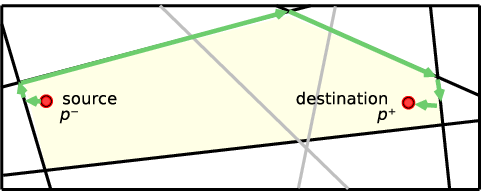}

\caption{The path marked by arrows
illustrates one of the two possible journeys around the cell
$\mathcal{C}(p^-,p^+)$ which start at source $p^-$ and end at
destination $p^+$.}\label{fig:construction}
\end{figure}

Any two specified points $p^-$ and $p^+$ in the plane will almost
surely not be hit by any of the lines of a given isotropic stationary
Poisson line process $\Pi$ and will therefore fail to be connected by
$\Pi$. Accordingly, we establish the convention that movement from
$p^-$ to $p^+$ occurs as follows: first, use the Poisson tessellation
to construct the cell $\mathcal{C}(p^-,p^+)$ containing $p^-$ and
$p^+$ which arises by deleting all Poisson lines which separate $p^-$
from $p^+$. Now, proceed from the source $p^-$ in exactly the opposite
direction to that of $p^+$ until one first encounters a Poisson line
[which will be part of the cell boundary $\partial\mathcal
{C}(p^-,p^+)$]. Then, continue along the line in one or the other
direction, clockwise or counterclockwise, proceeding along the boundary
of the cell $\mathcal{C}(p^-, p^+)$. Continue until one reaches the
ray extending from $p^-$ and through $p^+$. Then, proceed down this ray
to the destination $p^+$ (Figure~\ref{fig:construction}). Thus, we
consider \textit{near-geodesics}; routes based on \textit{semiperimeters}
of the cell $\mathcal{C}(p^-, p^+)$. These are to be considered in
contrast to \textit{network geodesics}, which always use the shortest
network path and can therefore be found only by solving a difficult
optimization problem.

%

Evidently, this is a conservative option for plumbing nodes into the
Poisson network produced by $\Pi$, suitable if we wish to produce\vadjust{\goodbreak}
upper bounds on connection lengths and adopted without further comment
in what follows. We suppose the choice of whether to travel clockwise
or counterclockwise around the cell $\mathcal{C}(p^-,p^+)$ (or,
equivalently, which semiperimeter to choose) is made at random and
equiprobably, independently for each pair of nodes $p^-$, $p^+$.
(As mentioned in Section~\ref{sec:conclusion}, interesting and hard
problems arise if the choice of route for a specific pair is influenced
by the flow in the entire network.)

In contrast to true network geodesic connections, these routes can be
viewed as outputs from an unsophisticated but direct \textit
{semiperimeter algorithm}: if one is on a Poisson line and encounters
another Poisson line, then one chooses (from the three onward paths)
that path which leads closest to the eventual destination without
separating source from destination. This focuses attention on the
\textit
{Poissonian city}, a region connected by routes based on a fixed
stationary isotropic Poisson line process and following the above
convention so as to ensure that the line process actually connects nodes.
Questions addressed in Section~\ref{sec:geodesics} of this paper,
filling in and extending the results announced in
\citet{Kendall-2009b}, include the following:
\begin{enumerate}
\item What can one say about the basic geometry of these routes?
Computations from \citeauthor{AldousKendall-2007} [(\citeyear
{AldousKendall-2007}), Theorem 3] yield
$\frac{4}{3}\log\dist(p^-,p^+)$ as asymptotic mean excess length as
$\dist(p^-,p^+)$ tends to infinity. This can be viewed as a
quantitative development of the announcement by
\citeauthor{Davidson-1974c} [(\citeyear{Davidson-1974c}), Theorem
5(ii)], but by how much does the traveled path deviate
laterally from the Euclidean connection, and at what point is that
lateral deviation greatest? (See Theorem~\ref{thm:lateral-displacement}.)
\item What is the order of random variation of the route lengths? (See
Theorem~\ref{thm:variance-bound}.)
\item What might be said about how actual network geodesics differ from
these routes? (This is discussed in Section~\ref{sec:true-geodesics}.)
\item Finally, might actual network geodesics produce substantially
smaller mean excess lengths? (Theorem~\ref{thm:poisson-lower-bound}
shows that the mean excess of true network geodesics is comparable to
the mean excess of semiperimeter paths.)
\end{enumerate}


\subsection{Traffic flow}\label{sec:traffic-flow}

Given a Poissonian city, it is natural to consider traffic flow.
Suppose, for example, that the city is represented by
$\ball(\origin,n)$, a disc centered at $\origin$ and of radius
$n$, connected by roads provided by a stationary isotropic Poisson
line process $\Pi$. Suppose that each pair of points $p^-$ and
$p^+$ in the disc generates a constant infinitesimal amount of
traffic $\dd p^-\,\dd p^+$, divided equally between each of the two
connecting routes generated according to the semiperimeter algorithm
described above. Suppose further that we condition on the event of a
Poisson line passing through the center $\origin$ of the disc. Let
%
%
\begin{equation}\label{eq:region}
\mathcal{D}_n =
\{
(p^-, p^+)\in\ball(\origin,n)^2 \dvtx p^-_1<p^+_1, \origin\in
\partial\mathcal{C}(p^-,p^+)
\}
\end{equation}
denote the region in $4$-space corresponding to $p^-$, $p^+$ in $\ball
(\origin,n)$ for which
$p^-$ lies to the left of $p^+$ (imposed by the inequality
$p^-_1<p^+_1$, where $p^-_1$, $p^+_1$ are the $x$-coordinates of source
and destination nodes) and
one of the two possible routes passes through $\origin$.


Questions addressed in Section~\ref{sec:flow} include:
\begin{enumerate}
\item What is the dependence on $n$ of the mean $\mathbb{E}[T_n]$ of
\begin{eqnarray*}
T_n&=&
\frac{1}{2}\iint\mathbb{I}_{[(p^-,p^+)\in\mathcal{D}_n]}\,\dd
p^-\,\dd p^+
\\
&=&
\frac{1}{2}\iint_{\ball(\origin,n)^2}\mathbb{I}_{[p^-_1<p^+_1,
\origin\in\partial\mathcal{C}(p^-,p^+)]}\,\dd p^-\,\dd p^+ ,
\end{eqnarray*}
the total amount of traffic
passing through the center $\origin$? This quantity scales as~$n^3$,
following from scaling arguments using basic stochastic geometry.
However, one can in fact compute the constant of asymptotic
proportionality (Theorem~\ref{thm:first-moment}).
\item Indeed, Aldous has asked whether the scaled flows $T_n/n^3$ have
a nondegenerate limiting distribution. (The answer is that they do: see
Theorem~\ref{thm:scaling-limit} and Corollary~\ref{cor:non-degenerate}.)
\item Does uniform integrability hold for the sequence of $T_n/n^3$ as
$n\to\infty$? If not, then there might exist a well behaved limiting
distribution, but the mean of $T_n/n^3$ might converge to a higher
value than that of the limit. Were this the case, it could be viewed as
a kind of stochastic congestion result.
(The results of Theorem~\ref{thm:first-moment} and Lemma \ref
{thm:scaling-limit-mean} indicate why uniform integrability does hold;
it is possible to push this further, as indicated in Section \ref
{sec:uniform-integrability}.)
\end{enumerate}

Section~\ref{sec:grid-structure} provides a comparison by giving an
overview of analogous results for flows in cities built on grids
(\textit
{Manhattan cities}).
The concluding Section~\ref{sec:conclusion} adds some further remarks
and mentions possible future research directions.

\subsection{Directory of results}\label{sec:directory}

Finally, we present a directory of the main results so as to assist the
reader in navigating around this paper.

In Section~\ref{sec:connections}, Theorem \ref
{thm:lateral-displacement} establishes statistical asymptotics for the
maximum lateral deviation of a near-geodesic from the corresponding
Euclidean path; asymptotically, the maximum will be located uniformly
along the path, with extent given by the radial part of a
four-dimensional Gaussian vector with variance which is quadratic in the
location of the maximum and which vanishes at the two endpoints. A
preliminary Lemma~\ref{lem:tau-asymptotics} then leads to Theorem \ref
{thm:variance-bound}, which produces an asymptotic bound $\frac
{20}{27}\log n$ for the variance of the excess of near-geodesic length
over Euclidean path length, for a near-geodesic started at a point and
going off to a line at distance $n$ from the point. Theorem \ref
{thm:poisson-lower-bound} generates an asymptotic lower bound $2(\log4
- \frac{5}{4})\log n$ on the mean excess for any path between two
fixed points separated by distance $n$. The numerical value
$0.27258872\ldots$ of the constant of proportionality here should be
compared with the corresponding constant $\frac{4}{3}=1.33\ldots$
for near-geodesics [\citet{AldousKendall-2007}].

In Section~\ref{sec:flow} the focus changes to flows in networks built
from Poisson line processes. After introducing the notion of a
``Poissonian city'' (based on a disc of radius $n$), Theorem \ref
{thm:first-moment} shows that the traffic flow $T_n$ through the center
has asymptotic mean $2n^3$. Corollary~\ref{cor:sharper-asymptotic}
reports a refinement of the detailed asymptotics; Theorem \ref
{thm:scaling-limit} (supported by Lemma~\ref{thm:scaling-limit-mean})
establishes the existence of a limiting distribution for
$T_n/n^3$---although the result only gives a geometric
characterization---and Corollary
\ref{cor:non-degenerate} confirms that this limit distribution is
nondegenerate.

\section{Making connections in the Poissonian city}\label{sec:geodesics}
Asymptotic arguments applied to formulas from stochastic geometry
indicate both the geometry of the routes provided by the
unsophisticated semiperimeter algorithm described above (including the
extent of random variation in length) and also ways in which they
differ from true network geodesics between source and destination
nodes. We begin by discussing the asymptotic distribution of the
location and extent of the maximum lateral displacement of a
semiperimeter route from the corresponding Euclidean geodesic.

\subsection{Maximum lateral displacement}\label{sec:maximum}

%

\begin{figure}

\includegraphics{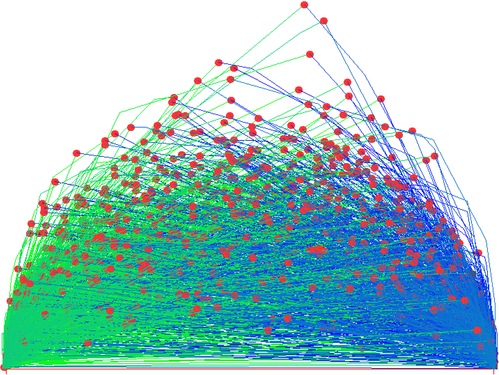}

\caption{A plot of $1000$ semiperimeters
of cells $\partial\mathcal{C}(p^-,p^+)$ based on a distance
$n=\dist(p^-$, $p^+)=1000$. The dots indicate the maximum lateral
displacements from the horizontal axis between source and destination.
The figure has been subjected to vertical exaggeration by a factor of
$\sqrt{n}/4$.}\label{fig:scale}
\end{figure}

Consider the height and location of the maximum lateral displacement of
one of the $\partial\mathcal{C}(p^-, p^+)$ semiperimeter routes from
a source $p^-$ to a destination $p^+$. Figure~\ref{fig:scale}
%
%
illustrates $1000$ realizations of such routes, with maxima marked by
discs, when source and destination are separated by distance $n=1000$.
Such simulations suggest
the existence of a limiting distribution under scaling for the extent
and location of the maximum lateral displacement, and stochastic
geometry arguments show that this is indeed the case.
\begin{theorem}\label{thm:lateral-displacement}
Consider two points $p^-=\origin=(0,0)$ and $p^+=(n,0)$ located along
the $x$-axis and also a path between these points based on $\partial
\mathcal{C}(p^-,p^+)\cap\{(x,y)\dvtx y\geq0\}$. Locate the maximum
lateral displacement of $\partial\mathcal{C}(p^-,p^+)\cap\{
(x,y)\dvtx y\geq0\}$ from the $x$-axis (and thus from the Euclidean
geodesic between $p^-$ and $p^+$) at $(n U_n,\sqrt{n} V_n)$. The joint
distribution of $(U_n, V_n)$ then has the following weak limit $(U,V)$
as $n\to\infty$:
\begin{enumerate}[(a)]
\item[(a)] the scaled location $U$ is uniformly distributed over
$[0,1]$;
\item[(b)] conditional on $U=u\in(0,1)$, the
scaled displacement $V$ is distributed as the length of a
four-dimensional Gaussian vector with variance \mbox{$2 u(1-u)$}.
\end{enumerate}
\end{theorem}
%
%
\begin{pf}
Let $p^-$, $p^+$ be located at $\origin$ and $(n,0)$ on the $x$-axis
and let the maximum displacement be located\vadjust{\goodbreak} at $(n U_n, \sqrt{n}V_n)$,
as in the statement of the theorem, so that
\begin{eqnarray*}
\sqrt{n}V_n &=& \max\{y \dvtx(x,y)\in\partial
\mathcal{C}(p^-,p^+)\},\\
\bigl(n U_n, \sqrt{n}V_n\bigr) &\in& \partial\mathcal{C}(p^-,p^+) .
\end{eqnarray*}
(Almost sure uniqueness of $U_n$ is a consequence of the fact that a
stationary isotropic Poisson line process almost surely contains no
horizontal lines.)

The proof is a variation on ideas in the proof of Theorem 3 in
\citet{AldousKendall-2007}.
Consider the point process formed by intersections of lines $\ell^-$,
$\ell^+$ from $\Pi$, subject to the following, additional, requirements:
\begin{enumerate}
\item
no further lines from $\Pi$ separate the intersection $\ell^-\cap
\ell^+$ from the segment of length $n$ formed between the pair of
points $p^-=\origin$, $p^+=(n,0)$;
\item
one of the intersecting lines $\ell^-$ has positive slope, the other
$\ell^+$ has negative slope and neither line intersects the segment
formed between the pair of points $p^-$, $p^+$.
\end{enumerate}
Topological arguments show that there must be just two points in this
point process, one above and one below the $x$-axis, and the point
above the $x$-axis must be located at $(n U_n,\sqrt{n}V_n)$. The
intensity of the point process in the upper half-plane is given by
%
%
\begin{eqnarray}
\rho(x,y) &=&
\tfrac{1}{4}
\bigl(\sin\alpha+\sin\beta-\sin(\alpha+\beta)\bigr)\nonumber\\[-8pt]\\[-8pt]
&&{}\times\exp\bigl(-\tfrac{1}{2}\bigl(\sqrt{x^2+y^2}+\sqrt
{(n-x)^2+y^2}-n\bigr)\bigr),
\nonumber
\end{eqnarray}
where $\alpha$, $\beta\in(0,\pi)$ are the interior angles at
$\origin$ and $(n,0)$ of the triangle formed by $(x,y)$ and these two points.
Here, the exponential factor is contributed by requirement
1 above, since Slivynak's theorem
can be used to show that the unit intensity stationary isotropic
Poisson property is preserved by conditioning on two lines from $\Pi$
intersecting at $(x,y)$ and then removing those two lines. Employing
the fact that the intensity of the point process formed by
intersections of lines from the unit intensity line process $\Pi$ is
$\frac{\pi}{2}$, requirement 2 can be shown
to lead to the factor
\begin{eqnarray*}
&&\frac{\pi}{2}\times\frac{1}{\pi}\int_0^\alpha\int_0^{\beta
}\frac{1}{2}\sin(\theta+\psi) \,\dd\psi\,\dd\theta\\
&&\qquad=
\frac{1}{4}
\bigl(\sin\alpha+\sin\beta-\sin(\alpha+\beta)\bigr) .
\end{eqnarray*}
%

It follows that $(U_n,V_n)$ has joint density on the upper half-plane
given asymptotically for large $n$ when $0<u<1$ and $v>0$ by
\begin{eqnarray*}
n^{3/2}\rho\bigl(nu,\sqrt{n}v\bigr)
&=&
\frac{n^{3/2}}{4}\bigl(\sin\alpha+ \sin\beta-
\sin(\alpha+\beta)\bigr)\\
&&{}\times\exp\Biggl(-\frac{n}{2}\Biggl(\sqrt{u^2+\frac{v^2}{n}}+\sqrt
{(1-u)^2+\frac{v^2}{n}}-1\Biggr)\Biggr)\\
&\sim&
\frac{n^{3/2}}{4}
\exp\biggl(-\frac{1}{4}\frac{v^2}{u(1-u)}\biggr)\\
&&{}\times
\bigl(\sin(\alpha)\bigl(1-\cos(\beta)\bigr) + \sin(\beta)\bigl(1 - \cos
(\alpha)\bigr)\bigr) .
\end{eqnarray*}
Converting the sines and cosines to expressions in $u$ and $v$, as
$n\to\infty$
\begin{eqnarray*}
&&\frac{n^{3/2}}{4}\rho\bigl(nu,\sqrt{n}v\bigr) \\
&&\qquad\sim
\frac{n^{3/2}}{4}\exp\biggl(-\frac{1}{4}\frac{v^2}{u(1-u)}\biggr)\\
&&\qquad\quad{}\times\biggl(\frac{{v}/{\sqrt{n}}}{\sqrt{u^2+{v^2}/{n}}}
\biggl(1-\frac{1-u}{\sqrt{(1-u)^2+{v^2}/{n}}}\biggr)\\
&&\qquad\quad\hspace*{17.8pt}{} + \frac{{v}/{\sqrt{n}}}{\sqrt
{(1-u)^2+{v^2}/{n}}}
\biggl(1-\frac{u}{\sqrt{u^2+{v^2}/{n}}}\biggr)\biggr)\\
&&\qquad\to
\frac{1}{8}\frac{v^3}{u^2(1-u)^2}\exp\biggl(-\frac{1}{4}\frac
{v^2}{u(1-u)}\biggr) .
\end{eqnarray*}
This can be identified as the joint density corresponding to the
limiting distribution of $U_n$ and $V_n$ given in the theorem; weak
convergence follows from Fatou's lemma.
\end{pf}
%

Simulation studies (from which Figure~\ref{fig:scale} was derived)
confirm these asymptotics.

\subsection{Random variation via growth processes}\label{sec:growth}
The direct stochastic geometry method is highly effective for computing
detailed asymptotics of mean-value quantities but leads to burdensome
calculations for second order quantities such as variances. We
therefore turn to an alternate approach based on random growth processes.
The maximum analyzed in Section~\ref{sec:maximum} occurs at the point
of intersection of the trajectories of two independent \textit{growth
processes}, together representing the lateral deviation of the path
from the Euclidean geodesic between source and destination nodes. One
growth process is viewed as starting from the source node $p^-$ and one
from the destination node $p^+$, tracing out the relevant semiperimeter
of $\mathcal{C}(p^-,p^+)$ by describing the height as a function of
arc length along the semiperimeter. The two processes $\{H^\pm_s\dvtx
s\geq0\}$ are given by heights $H^\pm_s$ above the $y$-axis at arc
length distance $s$ along the respective path from the originating node
($p^-$ for $H^+$, $p^+$ for $H^-$).
Let $\Theta^\pm_s$ be the angle made by the path with the $x$-axis,
where the angle measurement is oriented depending on the label ``$\pm
$'' so that $\Theta^\pm_0=\pi$ (since the path commences by setting
out in the opposite direction to that of its goal). Then, $\frac{\dd
}{\dd s}H^\pm_s=\sin\Theta^\pm_s$ [except for isolated points at
which the slope of $\partial\mathcal{C}(p^-, p^+)$ changes]. Changes
in $\Theta^\pm$ occur when the path is intercepted by a line from
$\Pi$ which also intersects that part of the $x$-axis with $p^\mp$
deleted which does not contain $p^\pm$. Indeed, starting at arc length
$s$ along the path from $p^\mp$, the angle $\Theta^\pm$ remains
constant for an Exponentially distributed\vspace*{1pt} length of rate $\frac
{1}{2}(1-\cos\Theta^\pm_s)$, after which\vspace*{1pt} the angle jumps to a lower
value $\Theta^\pm_{s-}+\Delta\Theta^\pm_s<\Theta^\pm_{s-}$. The
hitting properties of Poisson line processes can be used to show
%
%
\begin{equation}\label{eq:theta-jump}
\mathbb{P}[-\Delta\Theta^\pm_s\leq\phi| \Theta^\pm_{s-}=\theta
] = \frac{1-\cos\phi}{1-\cos\theta}\qquad
\mbox{where }0\leq\phi\leq\theta.
\end{equation}
(Jump processes such as $\Theta$ are taken to be c\`adl\`ag so that
$\lim_{s\searrow t}\Theta_s=\Theta_t$, and we write $\Theta_-$ for
the process of left limits.)
We suppose the growth processes $H^\pm$ evolve for all time according
to the dynamics described above.
Let $X^\pm_s$ be the distance from $p^\mp$ when resolved along the
axis from source $p^\mp$ to destination $p^\pm$,
so $\frac{\dd}{\dd s}X^\pm_s=\sec\Theta^\pm_s$ (except for
isolated points of discontinuity of $\Theta$).
%
%
\begin{figure}

\includegraphics{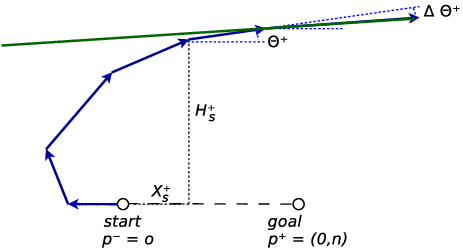}

\caption{Illustration of $H^+$
construction. The growth process $H^+_s$ tracks height as a function
of arc length $s$. The angle of slope $\Theta^+_s$ is an auxiliary
process governed by a Poisson stochastic differential equation,
jumping when the path is intercepted by a line from $\Pi$ which also
intersects the negative
$x$-axis.}\label{fig:growth}
\end{figure}
Figure~\ref{fig:growth} illustrates the $H^+$ construction.

\subsubsection*{Analysis of growth process using L\'evy processes}%

For the sake of clarity of exposition, we now drop the superscripted
$\pm$.

It is convenient to apply a random time change $t=s-X_s$ using the
excess of arc length over distance traveled toward the goal; in the new
time scale, the angle $\Theta$ changes according to a Poisson process
of incidents of rate $\frac{1}{2}$ while $\frac{\dd H}{\dd t}=\sin
\Theta/(1-\cos\Theta)$, $\frac{\dd X}{\dd t}=\cos\Theta/(1-\cos
\Theta)$. This gives a stochastic differential equation for $H$ and\vadjust{\goodbreak}
$X$, driven indirectly by a half-unit rate Poisson counting process $N$
via the auxiliary process $\Theta$:
%
%
\begin{eqnarray}\label{eq:Poisson-sde}
\dd H &=& \frac{\sin\Theta}{1-\cos\Theta}\,\dd t ;\nonumber
\\
\dd X &=& \frac{\cos\Theta}{1-\cos\Theta}\,\dd t ;
\\
\dd\Theta&=& \Delta\Theta\,\dd N .\nonumber
\end{eqnarray}
The distribution of the jump $\Delta\Theta$ is given by (\ref
{eq:theta-jump}), and the jump $\Delta\Theta$ at a jump of $N$ is
conditionally independent of the past given $\Theta_-$ at that time.
Note that~(\ref{eq:Poisson-sde}) can be viewed as driven by a \textit
{marked} Poisson process, obtained by marking the incidents of $N$ by
the jumps of $\Theta$, with mark distribution (conditional on the left
limit $\Theta_-$) given by~(\ref{eq:theta-jump}).
Using\vspace*{1pt} this terminology and approach, we will now show that
the excess
$\sigma(n)=\inf\{t\dvtx X_t\geq n\}=\inf\{t\dvtx\int_0^t\frac{\cos\Theta
_u}{1-\cos\Theta_u}\,\dd u\geq n\}$ at given distance $X=n$ has standard
deviation asymptotically proportional to $\sqrt{\log n}$ for large $n$.

%

%
%
\begin{figure}

\includegraphics{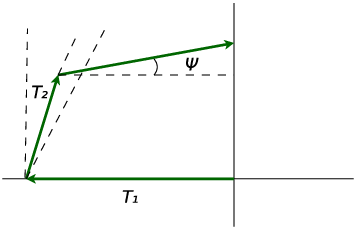}

\caption{Illustration of initial segment
of a path used to provide an upper bound on the excess acquired before
$X=0$. This initial segment is made up of three line segments, of
lengths $T_1$, $T_2$, and finally a length bounded above by $T_1
\sec U$. Here, $T_1$, $T_2$, $U$ are independent with
distributions as given in the text.}\label{fig:excess}
\end{figure}

To establish this result it
is simplest to consider the growth process begun with $X_0=0$ and
$\Theta_0$ lying in the range $(0, \pi/2]$. We must therefore control
the amount of excess required to achieve this over the initial segment
of the path. Since $X=0$ at both ends of this segment, this excess can
be bounded above by the length of the initial segment of the path
indicated in Figure~\ref{fig:excess}. This path uses the following
directions until it first hits the $y$-axis: it first runs along the
negative $x$-axis until it encounters a line of angle between $\frac
{\pi}{3}$ and $\frac{\pi}{2}$; it then moves upward along this line
until it encounters a line of angle between $0$ and $\frac{\pi}{3}$,
then along that line. The first segment is of length $T_1$,
Exponentially distributed of rate $\frac{1}{4}$. The second segment
is of length stochastically bounded above by $T_2$, Exponentially
distributed of rate $\frac{1}{4}$. The final segment is of length
stochastically bounded above by $T_1\sec U$, where $U$ has density
$\frac{2}{\sqrt{3}}\cos u$ for $0<u<\frac{\pi}{3}$ (this uses a
stochastic monotonicity argument applied to the conditional
distribution of the angle of the third line segment given the angle of
the second line segment). We may take $T_1$, $T_2$, $U$ to be independent.

The quantity $T_1+T_2+T_1\sec U$ is an upper bound on the path length
of the initial segment and has finite mean given by $8(1+\frac{\pi
}{3\sqrt{3}})$ and finite\vspace*{-1pt} second moment given by $32 (3+\frac
{2}{\sqrt{3}}( \pi+ \log(2+\sqrt{3})))$.
It follows that the contribution of the actual initial segment to the
mean and variance of the excess is bounded and may be ignored if we can
establish logarithmic increase in variance of the remainder.
Accordingly, we may suppose our growth process begins with $X_0=0$ and
with $\Theta_0$ lying between $0$ and $\frac{\pi}{2}$, distributed
according to the density $\cos\theta$ for $0<\theta<\frac{\pi
}{2}$. [This is because the growth process will intersect the positive
$y$-axis in the first intercept which makes an angle with the
horizontal in the range $(0,\frac{\pi}{2})$.]

We now address the question of the variance of $\sigma(n)$ (based on
$X_0=0$) in three stages. First, we use trigonometry and probabilistic
coupling to relate the negative log angle $-\log\Theta$ to a Le\'vy
subordinator $\xi$. We then state and prove a lemma on analogous mean
and variance asymptotics for $\tau(n)=\inf\{t\dvtx\int_0^t\exp(2\xi
_u)\,\dd u\geq n\}$. Finally, we state and prove a theorem which uses
approximations and coupling to establish the required mean and variance
asymptotics for $\sigma(n)$.

For the first stage, note that if $0\leq\Theta\leq\frac{\pi}{2}$,
then, by trigonometry and calculus, we can convert the first two
partial sums
of the Laurent series for $\frac{\cos\Theta}{1-\cos\Theta}$ into
upper and lower bounds over $[0,\frac{\pi}{2}]$:
%
%
\begin{equation}\label{eq:trig-bound-1}
\frac{2}{\Theta^2} -\frac{5}{6} \leq\frac{\cos\Theta
}{1-\cos\Theta} \leq\frac{2}{\Theta^2} .
\end{equation}
The application of the distributional information in (\ref
{eq:theta-jump}) to a jump $\Delta\Theta=\Theta-\Theta_-$ produces
a unit Exponential random variable:
%
%
\begin{equation}\label{eq:mark-construction}
\mathcal{J} = - \log\biggl(\frac{1-\cos(-\Delta\Theta
)}{1-\cos\Theta_-}\biggr) .
\end{equation}
We can thus mark independently each jump of the Poisson process $N$
using the unit Exponential mark distribution of $\mathcal{J}$. The
mark $\mathcal{J}$ can be used, together with~$\Theta_-$, to
reconstruct the actual jump of $\Theta$, using~(\ref{eq:mark-construction}).
Bearing in mind that $\Delta\Theta<0$, we may write
\[
\mathcal{J} = f(\log\Theta_-)-f(\log(-\Delta\Theta)),
\]
where $f(x)=\log(1-\cos e^{x})$. Calculus shows that $f^\prime>0$ and
$f^{\prime\prime}<0$ over the range $(-\infty,\log\frac{\pi
}{2})$. Hence, if $x \leq\log\Theta_-\leq\log\frac{\pi}{2}$, then
\[
\frac{\pi}{2} \leq
f^\prime(\log\Theta_-) \leq f^\prime(x)
=
\frac{e^x \sin e^x}{1 - \cos e^x}
< 2 = \lim_{x\to-\infty}\frac{e^x\sin e^x}{1-\cos e^x} .
\]
However, $0<-\Delta\Theta<\Theta_-$ and $\mathcal{J}=\int_{\log
(-\Delta\Theta)}^{\log\Theta_-}f^\prime(x)\,\dd x$; therefore,
\[
f^\prime(\log\Theta_-)\bigl(\log\Theta_- - \log(-\Delta\Theta
)\bigr)
\leq
\mathcal{J}
\leq
2 \bigl(\log\Theta_- - \log(-\Delta\Theta)\bigr) .
\]
Hence, for $0<-\Delta\Theta<\Theta_-\leq\frac{\pi}{2}$,
%
%
\begin{eqnarray}\label{eq:Phi-inequality}
-\log\biggl(1-\exp\biggl(-\frac{2}{\pi}\mathcal{J}\biggr)\biggr) &\leq&
-\log\biggl(1-\exp\biggl(-\frac{\mathcal{J}}{f^\prime(\log\Theta
_-)}\biggr)\biggr)\nonumber\\
&\leq&
-\log\biggl(\frac{\Theta_-+\Delta\Theta}{\Theta_-}\biggr) =
-\Delta\log\Theta\\
&\leq&
-\log\biggl(1-\exp\biggl(-\frac{1}{2}\mathcal{J}\biggr)\biggr) .
\nonumber
\end{eqnarray}
%
For future reference, note that $-\log(1-\exp(-\frac
{1}{2}\mathcal{J}))$ is distributed as the maximum $T^\prime
\vee T^{\prime\prime}$ of two independent unit-mean Exponential
random variables $T^\prime$ and $T^{\prime\prime}$, and, in
particular, it has mean $\frac{3}{2}$ and variance $\frac{5}{4}$.
On the other hand
$-\log(1-\exp(-\frac{2}{\pi}\mathcal{J}))$ has
probability density $\frac{\pi}{2}(1-e^{-x})^{{\pi
}/{2}-1}e^{-x}$ for $x>0$, and the square of its negative exponential
$(1-\exp(-\frac{2}{\pi}\mathcal{J}))^2$
(which will play a role later on) has probability density
%
%
\begin{equation}\label{eq:lowerbound-density}
\frac{\pi}{4}\bigl(1-\sqrt{x}\bigr)^{{\pi}/{2}-1}\frac{1}{\sqrt{x}}
\qquad\mbox{for }0<x<1 .
\end{equation}

Thus, a coupling construction indicated by the inequalities of (\ref
{eq:Phi-inequality}) permits approximation of the negative logarithm of
the angle process by
%
$\eta\leq-\log(\Theta/\Theta_0)\leq\xi$.
Here, $\eta$, $\xi$ are nondecreasing pure-jump L\'evy processes
(hence subordinators) which have jumps at the same times as those of
$-\log\Theta$
(namely, at\vadjust{\goodbreak} incidents of the Poisson counting process $N$ of intensity
$\frac{1}{2}$),
but with jump distributions given by the distributions of $-\log
(1-\exp(-\frac{2}{\pi}\mathcal{J}))$ and $-\log
(1-\exp(-\mathcal{J}/2))$, respectively. For future reference, note
also that the Laplace exponent $\Phi(q)=-\frac{1}{t}\log\mathbb
{E}[e^{-q \xi_t}]$ of $\xi$ can be computed as $\Phi(q)=\frac
{q(3+q)}{2(1+q)(2+q)}$ for $q>-1$, while $M_s=\xi_s -\frac{3}{4}s$
defines a martingale, as does $M^2_s-\frac{5}{8} s$. In particular,
$M$ is an $L^2$-martingale.

The discrepancy between coupled jumps of $\xi$ and $-\log\Theta$ can
be controlled by
%
%
\begin{equation}\label{eq:coupling-discrepancy}
0 \leq
\Delta\xi- \biggl(-\log\frac{\Theta_-+\Delta\Theta}{\Theta
_-}\biggr)
\leq
\log\biggl(\frac{1-\exp(-\mathcal{J}/f^\prime(\log\Theta
_-))}{1-\exp(-{\mathcal{J}}/{2})}\biggr) .
\end{equation}
Now, if $0\leq a\leq\frac{\pi}{2}$, then we can use the
inequalities $\sin\frac{a}{2}\leq\frac{a}{2}$, $\cos^2\frac
{a}{2}\geq\frac{1}{2}$ and the convexity of $\tan^2\frac{a}{2}$
over this range to establish the simple bound
\[
\frac{1}{f^\prime(\log a)} = \frac{1-\cos a}{a\sin a}
\leq\frac{1-\cos a}{\sin^2 a} =
\frac{1}{2}+\frac{\tan^2({a}/{2})}{2}
\leq
\frac{1}{2}+\frac{a^2}{4}
\]
%
and so (using $\Theta\leq\frac{\pi}{2}$ to apply the above inequality)
\begin{eqnarray*}
\log\biggl(\frac{1-\exp(-\mathcal{J}/f^\prime(\log\Theta
_-))}{1-\exp(-{\mathcal{J}}/{2})}\biggr)
&\leq&
\log\biggl(1+\frac{1-\exp(-({\Theta_-^2}/{2})
({\mathcal{J}}/{2}))}{\exp({\mathcal{J}}/{2}
)-1}\biggr)\\
&\leq&
\log\biggl(1+\frac{{\mathcal{J}}/{2}}{\exp(
{\mathcal{J}}/{2})-1}\frac{\Theta_-^2}{2}\biggr)
\\
&\leq&
\frac{{\mathcal{J}}/{2}}{\exp({\mathcal
{J}}/{2})-1} \frac{\Theta_-^2}{2}
\leq
\frac{\Theta_-^2}{2}\\
&\leq&
\frac{1}{2}\exp(-2\eta_-).
\end{eqnarray*}

Applying this upper bound on the jumps, it follows that we can control
the total discrepancy between $\xi$ and $-\log\Theta$ by
%
%
\begin{equation}\label{eq:final-discrepancy}
0
\leq
\xi_t - \bigl(- \log(\Theta_t/\Theta_0)\bigr)
\leq
\frac{1}{2} \mathop{\sum_{w\leq t:}}_{\Delta N_w>0} \exp
(-2\eta_{w-})
.
\end{equation}
However, the right-hand side is increasing in $t$ and has a limit which
can be expressed in terms of a simple perpetuity. Indeed,
\[
\sum_{w \dvtx\Delta N_w>0}
\exp(-2\eta_{w-})
= U_t \leq U_\infty=
\Biggl(1 + \sum_{k=1}^\infty\prod_{m=1}^k\exp(-2\Delta_m\eta
)\Biggr) ,
\]
where $\Delta_m\eta$ is the $m$th jump of $\eta$.
A classical calculation gives the first and second moments of the
perpetuity final value $U_\infty$ in terms of first and second moments
of $\exp(-2\Delta_m\eta)$
[see, e.g., \citet{Vervaat-1979}, Theorem~5.1]. However, we will\vadjust{\goodbreak}
require control of exponential moments $\mathbb{E}[\exp(z U_\infty)]$
for positive~$z$. \citet{AlsmeyerIksanovRosler-2009} and
\citet{Kellerer-1992c} give results for the general case, but in
our particular case the perpetuity multiplier $\exp(-2\Delta
_m\eta)$ is positive and bounded above by $1$; moreover, its
probability density~(\ref{eq:lowerbound-density}) is bounded above
near $1$, and it thus follows from monotonicity and the methods of
\citeauthor{GoldieGrubel-1996} [(\citeyear{GoldieGrubel-1996}), Theorem
3.1] that $\mathbb{E}[\exp(z
U_\infty)]<\infty$ for all positive $z$. [In fact, the upper bound
requirement on the probability density can be replaced by comparability
with a Beta density; see \citet{HitczenkoWesolowski-2010}, Section
4.]

We can now improve~(\ref{eq:trig-bound-1}) to provide bounds in terms
of L\'evy subordinators:
%
%
\begin{eqnarray}\label{eq:trig-bound-2}
\frac{2}{\Theta_0^2}\exp(2 \xi-U_\infty) -\frac{5}{6}
&\leq&
\frac{2}{\Theta^2} -\frac{5}{6} \leq
\frac{\cos\Theta}{1-\cos\Theta} \nonumber\\[-8pt]\\[-8pt]
&\leq&\frac{2}{\Theta^2}
\leq\frac{2}{\Theta_0^2}\exp(2 \xi) .
\nonumber
\end{eqnarray}
Accordingly, we first establish the following lemma, which delivers the
desired results for $\int\exp(2 \xi_s)\,\dd s$ rather than $X$.
\begin{lem}\label{lem:tau-asymptotics}
Define $\tau(n)$ in terms of $\int\exp(2 \xi_s)\,\dd s$ by
\[
n = \int_0^{\tau(n)} \exp(2 \xi_s)\,\dd s .
\]
Then,
%
%
\begin{equation}\label{eq:tau-representation}
\tau(n) = \frac{2}{3}\biggl(\log n - 2 M_{\tau(n)} + \log
\biggl(\exp\biggl(\frac{2\xi_{\tau(n)}}{n}\biggr)\biggr)\biggr)
\end{equation}
with the following asymptotics for mean and variance as $n\to\infty$:
%
%
\begin{eqnarray}
\label{eq:asymptmean}
\mathbb{E}[\tau(n)] &=& \tfrac{2}{3}\log n + O(1) ;\\
\label{eq:asymptvar}
\operatorname{Var}[\tau(n)] &=& \tfrac{20}{27}\log n + O\bigl(\sqrt{\log
n}\bigr) .
\end{eqnarray}
\end{lem}
\begin{pf}
First note the following trivial bound for $\tau(n)$ which establishes
the finiteness of $\mathbb{E}[\tau(n)]$:
\[
n = \int_0^{\tau(n)} \exp(2 \xi_s)\,\dd s
\geq\tau(n) .
\]

The representation~(\ref{eq:tau-representation}) was motivated by
heuristic time-reversal arguments and is essentially tautologous: write
\[
\exp\biggl(2M_{\tau(n)} + \frac{3}{2}\tau(n)\biggr) =
\exp\bigl(2\xi_{\tau(n)}\bigr) =
n \times\frac{\exp(2\xi_{\tau(n)})}{n} ,
\]
take logs and reexpress in terms of $\tau(n)$. Taking expectations, we
then obtain
%
%
\begin{equation}\label{eq:meantime}
\mathbb{E}[\tau(n)] = \frac{2}{3}\biggl(\log n + \mathbb{E}\biggl[\log
\biggl(\frac
{\exp(2\xi_{\tau(n)})}{n} \biggr)\biggr]\biggr) ;\vadjust{\goodbreak}
\end{equation}
the expectation $\mathbb{E}[M_{\tau(n)}]$ vanishes because $\tau(n)$
has finite expectation and so the stopped martingale $M_{s \wedge\tau
(n)}$ is $L^2$-bounded
by $\mathbb{E}[M_{\tau(n)}^2]\leq\frac{5}{8}\mathbb{E}[\tau(n)]$.

Consider the variance of $\tau(n)-\frac{2}{3} \log(\frac
{\exp(2\xi_{\tau(n)})}{n})$. Using $L^2$-martingale theory we find
that
%
%
\begin{equation}\label{eq:variancetime}
\operatorname{Var}\biggl[\tau(n)-\frac{2}{3} \log\biggl(\frac{\exp(2\xi
_{\tau(n)})}{n}\biggr)\biggr] = 4\biggl(\frac{2}{3}\biggr)^2
\operatorname{Var}\bigl[M_{\tau
(n)}\bigr] =
\frac{16}{9} \times\frac{5}{8} \mathbb{E}[\tau(n)] .\hspace*{-35pt}
\end{equation}
A uniform bound on the second moment of $\log(\frac{\exp(2\xi
_{\tau(n)})}{n})$ will therefore
permit us to deduce~(\ref{eq:asymptmean}) and~(\ref{eq:asymptvar})
from~(\ref{eq:meantime}) and~(\ref{eq:variancetime}), and so complete
the proof.

To this end, note that $Z_n=\exp(2\xi_{\tau(n)})$ constitutes a
Lamperti transformation [\citet{Lamperti-1972}] of the subordinator
$2\xi$ and therefore defines a self-similar Markov process $Z$. Using
the scaling property for $Z$ and maximizing $z (\log z)^2$ for $0<z<1$,
we obtain
%
%
\begin{eqnarray}\label{eq:low-bit}
\mathbb{E}\biggl[\biggl(\log\frac{\exp(2\xi_{\tau(n)})}{n}\biggr)^2 ;
\frac{\exp
(2\xi_{\tau(n)})}{n}<1\biggr]
&=&
\mathbb{E}\biggl[(\log Z_1)^2 ; Z_1<1 \Big| Z_0=\frac{1}{n}\biggr
]\hspace*{-30pt}\nonumber\\[-8pt]\\[-8pt]
&\leq&
4 e^{-2} \mathbb{E}\biggl[Z_1^{-1} \Big| Z_0=\frac{1}{n}\biggr] .
\nonumber
\end{eqnarray}
\citeauthor{BertoinYor-2005} [(\citeyear{BertoinYor-2005}), formula
(20), drawing on \citet{BertoinYor-2001}] can now be applied,
together with the calculation of the Laplace exponent $\Phi(q)$ of
$\xi$, to show that for $p>0$,
%
%
\begin{equation}\label{eq:higher-moments}
\mathbb{E}\biggl[Z_1^{-p} \Big| Z_0=\frac{1}{n}\biggr]
=
\frac{2p}{2p+1}\biggl(n^p e^{-n^p/2}+\biggl(\frac{n^p}{2}
\biggr)^{1-p}\int_0^{n^p/2} v^{p-1} e^{-v}\,\dd v\biggr) .\hspace*{-35pt}
\end{equation}
In particular, note that each positive integral moment of $Z_1^{-1}$ is
bounded, and, especially,
\[
\mathbb{E}\biggl[Z_1^{-1} \Big| Z_0=\frac{1}{n}\biggr] = \frac
{2}{3}\bigl(1+(n-1) e^{-n/2}\bigr) \leq\frac
{2}{3}(1+2e^{-{3}/{2}}) .
\]
Thus,~(\ref{eq:low-bit}) is bounded above.

To bound $\mathbb{E}[(\log\frac{\exp(2\xi_{\tau(n)})}{n} )^2 ;
\frac{\exp(2\xi_{\tau(n)})}{n}\geq1]$ for $n>0$, note that
$\tau(n)\leq\tau^\prime_n+1$, where $\tau^\prime_n=\inf\{t\dvtx2\xi
_t\geq\log n\}$. Hence,
\begin{eqnarray*}
&&\mathbb{E}\biggl[\biggl(\log\frac{\exp(2\xi_{\tau(n)})}{n}\biggr)^2 ;
\frac{\exp
(2\xi_{\tau(n)})}{n}\geq1\biggr] \\
&&\qquad\leq\mathbb{E}\biggl[(2\xi_{\tau^\prime_n+1}-\log n )^2 ; \frac
{\exp
(2\xi_{\tau(n)})}{n}\geq1\biggr]\\
&&\qquad\leq
\mathbb{E}[(2\xi_{\tau^\prime_n+1}-\log n)^2] .
\end{eqnarray*}
Now, $2\xi_{\tau^\prime_n+1}-\log n$ is the independent sum of a
summand of distribution $2\xi_1$ and a summand which is the overshoot\vadjust{\goodbreak}
$2\xi_{\tau^\prime_n}-\log n$. Since the jump distribution of $2\xi
$ is $2(T^\prime\vee T^{\prime\prime})$ (for two independent unit
Exponential random variables $T^\prime$ and~$T^{\prime\prime}$), it
follows by the memoryless property of the Exponential distribution that
the overshoot is some mixture of the distributions of $2 T^\prime$ and
$2(T^\prime\vee T^{\prime\prime})$, hence stochastically dominated
by $2(T^\prime\vee T^{\prime\prime})$. Consequently, it follows that
$\mathbb{E}[(2\xi_{\tau^\prime_n+1}-\log n)^2]\leq\frac
{111}{4}$.

Thus, $\mathbb{E}[(\log\frac{\exp(2\xi_{\tau(n)})}{n} )^2]\leq
\frac{2}{3}(1+2e^{-{3}/{2}})+\frac{111}{4}$, and
hence the lemma is proved.
\end{pf}

Note that upper bounds for $\mathbb{E}[(\log\frac{\exp(2\xi_{\tau
(n)})}{n})^2]$ can also be obtained using the techniques
of \citet{BertoinYor-2002}, but these bounds diverge to infinity
with $n$.

We are now able to deal with the asymptotic behaviors of the mean and
variance of $\sigma(n)$ as follows.
\begin{theorem}\label{thm:variance-bound}
With $\sigma(n)$ defined as above so that $n = X_{\sigma(n)}$ for $n>0$,
for large $n$,
%
%
\begin{eqnarray}
\mathbb{E}[\sigma(n)] &=& \tfrac{2}{3}\log n + O(1) , \\
\operatorname{Var}[\sigma(n)] &=& \tfrac{20}{27}\log n + O\bigl(\sqrt
{\log
n}\bigr) .
\end{eqnarray}
\end{theorem}
\begin{pf}
Lemma~\ref{lem:tau-asymptotics} provides partial control on the
asymptotic mean and variance of $\sigma(n)$ via (\ref
{eq:trig-bound-2}) since $X_t\leq\frac{2}{\Theta^2_0}\int^t_0 \exp
(2\xi_s)\,\dd s$, and so $\tau(\frac{\Theta^2_0}{2}n)\leq\sigma
(n)$. By the representation~(\ref{eq:tau-representation}) and
integration [since $\Theta_0$ has density $\cos\theta$ over
$(0,\frac{\pi}{2})$], we find that
%
%
\begin{eqnarray}\label{eq:upper-bound-asymptotics}\quad
&\displaystyle\mathbb{E}\biggl[\tau\biggl(\frac{\Theta^2_0}{2}n\biggr
)\biggr] = \frac
{2}{3}\log n +\frac{2}{3}\mathbb{E}\biggl[\log\biggl(\frac{\Theta
_0^2}{2}\biggr)\biggr]
+ O(1)
= \frac{2}{3}\log n + O(1) ,&\nonumber\\[-8pt]\\[-8pt]
&\displaystyle\operatorname{Var}\biggl[\tau\biggl(\frac{\Theta
^2_0}{2}n\biggr) -\frac{2}{3}\log
\biggl(\frac{\exp(2\xi_{\tau(\Theta^2_0 n)})}{\Theta^2_0 n}\biggr)\biggr]
= \frac{20}{27}\log n + O(1) .&\nonumber
\end{eqnarray}
The second moment of $\log(\frac{\exp(2\xi_{\tau(\Theta^2_0
n)})}{\Theta^2_0 n})$ being bounded,
it follows that
\[
\operatorname{Var}\biggl[\tau\biggl(\frac{\Theta^2_0}{2}n\biggr)\biggr]
=
\frac{20}{27}\log n + O\bigl(\sqrt{\log n}\bigr).
\]

We now consider the lower bounds from~(\ref{eq:trig-bound-2}),
delivering an upper bound for $\sigma(n)$
via
\[
\int_0^t \biggl(\frac{2}{\Theta_0^2}\exp(2 \xi_s - U_s) - \frac
{5}{6}\biggr) \,\dd s \leq
\int_0^t \biggl(\frac{2}{\Theta_s^2} - \frac{5}{6}\biggr) \,\dd s
\leq
X_t .
\]
First, observe what happens after the stopping time given by
\[
\kappa=
\biggl(1+\frac{5}{12}\Theta_0^2\biggr)\times\inf\biggl\{t \dvtx\frac
{2}{\Theta
_t^2}\geq\frac{2}{\Theta_0^2}+\frac{5}{6}\biggr\}.\vadjust{\goodbreak}
\]
We find that
\begin{eqnarray*}
\frac{2}{\Theta_{s+\kappa}^2}-\frac{5}{6} &\geq&
\frac{2}{\Theta_{0}^2} \frac{\Theta_{\kappa}^2}{\Theta_{s+\kappa
}^2} + \frac{5}{6} \frac{\Theta_{\kappa}^2}{\Theta_{s+\kappa
}^2}-\frac{5}{6}\\
&\geq&
\frac{2}{\Theta_{0}^2} \frac{\Theta_{\kappa}^2}{\Theta_{s+\kappa
}^2}
\end{eqnarray*}
and, moreover,
\[
\int_0^\kappa\biggl(\frac{2}{\Theta_s^2} - \frac{5}{6}\biggr) \,\dd
s \geq
-\frac{5}{6} \times\frac{\kappa}{1+({5}/{12})\Theta_0^2}
+ \frac{({5}/{12})\Theta_0^2}{1+({5}/{12})\Theta_0^2} \times
\frac{2\kappa}{\Theta_{0}^2}
= 0 .
\]
This will allow us to disregard the effects of the $-\frac{5}{6}$
term, so long as we can bound the second moment of $\kappa$. To this
end, introduce $\zeta$, a L\'evy subordinator jumping at the incidents
of the underlying Poisson process $N$, with jumps coupled to those of
the other processes so that $\zeta\leq\eta\leq- \log(\Theta
/\Theta_0)$, and with
\[
\Delta\zeta= \mathbb{I}_{[-\log(1-\exp(-({2}/{\pi})\mathcal
{J}))>1]} .
\]
Since these jumps are all of sizes $0$ or $1$, we can view $\zeta$ as
a Poisson counting process run at rate $\nu<\frac{1}{2}$. Moreover,
since $0<\Theta_0\leq\frac{\pi}{2}$, we have
\begin{eqnarray*}
\kappa&\leq&
\biggl(1+\frac{5}{12}\Theta_0^2\biggr)\inf\biggl\{t \dvtx\frac
{2}{\Theta_0^2}\exp(2\zeta_t)\geq\frac{2}{\Theta_0^2}+\frac
{5}{6}\biggr\}\\
&\leq&
\biggl(1+\frac{5\pi^2}{48}\biggr)\inf\biggl\{t \dvtx\zeta_t\geq\frac
{1}{2}\log\biggl(1+\frac{5\pi^2}{48}\biggr)\biggr\},
\end{eqnarray*}
and the boundedness of the second moment of the right-hand side of
these inequalities follows directly by comparison with a Gamma random
variable, derived from basic Poisson process properties. Hence,
$\mathbb{E}[\kappa^2]<\infty$.

So, consider $\widetilde{\xi}_s=\xi_{\kappa+s}-\xi_\kappa$ and
$\widetilde{\tau}(n)=\inf\{t\dvtx\int^t_0 \exp(2\widetilde{\xi
}_s)\,\dd
s=n\}$. We will bound $\sigma(n)$ above by a stopping time $\kappa+
\widetilde{\tau}(\frac{\Theta^2_0}{2}n)+\rho$, where $\rho$ is
chosen to compensate for the undershoot of $n$ at time $t=\kappa+
\widetilde{\tau}(\frac{\Theta^2_0}{2}n)$ caused by the $U$
contribution in
\[
\int_0^t \biggl(\frac{2}{\Theta_0^2}\exp(2 \xi_s - U_s) - \frac
{5}{6}\biggr) \,\dd s \leq
X_t .
\]
We find
\begin{eqnarray*}
&&\int_0^{\kappa+\widetilde{\tau}(({\Theta^2_0}/{2})n)+\rho}
\biggl(\frac{2}{\Theta^2}-\frac{5}{6}\biggr)\,\dd s\\
&&\qquad\geq
\int_0^{\widetilde{\tau}(({\Theta^2_0}/{2})n)+\rho}
\frac{2}{\Theta^2_0}\exp(2\widetilde{\xi}_s-U_{\kappa+s})\,\dd s\\
&&\qquad\geq
\exp\bigl(-U_{\kappa+\widetilde{\tau}(({\Theta
^2_0}/{2})n)}\bigr)
\biggl(n
+
\frac{2}{\Theta^2_0}\rho\exp\bigl(2\widetilde{\xi}_{\widetilde
{\tau}(({\Theta^2_0}/{2})n)}\bigr)\biggr) .
\end{eqnarray*}
Thus, we can choose $\rho$ to compensate for the undershoot so that
\[
\exp\bigl(-U_{\kappa+\widetilde{\tau}(({\Theta
^2_0}/{2})n)}\bigr) \biggl(n
+
\frac{2}{\Theta^2_0}\rho\exp\bigl(2\widetilde{\xi}_{\widetilde
{\tau}(({\Theta^2_0}/{2})n)}\bigr)\biggr)
=
n ;
\]
this is fulfilled by the choice
\[
\rho=
\bigl(\exp\bigl(U_{\kappa+\widetilde{\tau}(({\Theta
^2_0}/{2})n)}\bigr) - 1 \bigr)
\times\biggl(\frac{\exp(2\widetilde{\xi}_{\widetilde{\tau
}(({\Theta^2_0}/{2})n)})}{({\Theta^2_0}/{2})n}\biggr)^{-1} .
\]
Now, the first factor is bounded above by $\exp(U_\infty)$, and we
have already noted that $\mathbb{E}[\exp(z U_\infty)]<\infty$ for all
$z>0$ as a consequence of perpetuity theory.
The second factor is distributionally a randomization over $m$ of
$(\exp(2\xi_{\tau(m)})/m)^{-1}$,
and we have already noted that equation~(\ref{eq:higher-moments}), the
Lamperti transformation and the results of Bertoin and Yor allow us to
bound $\mathbb{E}[(\exp(2\xi_{\tau(m)})/m)^{-p}]$ uniformly in $m$ for
any fixed $p\geq1$.
We may thus deduce $\mathbb{E}[\rho^2]<\infty$ as a consequence of the
Cauchy--Schwarz inequality.

Accordingly, we find that
\[
\sigma(n) \leq\kappa+ \widetilde{\tau}\biggl(\frac{\Theta
^2_0}{2}n\biggr)+\rho
\]
bounds $\sigma(n)$ above by a stopping time with expectation $\frac
{2}{3}\log n+O(1)$; moreover, we may apply the representation (\ref
{eq:tau-representation}) to deduce that,
up to terms whose second moments are $O(1)$,
\[
\widetilde{\tau}\biggl(\frac{\Theta_0^2}{2}n\biggr) - \tau\biggl(\frac
{\Theta_0^2}{2}n\biggr)
\approx
2M_{\kappa+\widetilde{\tau}(({\Theta_0^2}/{2})n)} - 2M_{\tau
(({\Theta_0^2}/{2})n)},
\]
which itself must be of uniformly bounded second moment:
\begin{eqnarray*}
&&\mathbb{E}\bigl[\bigl( 2M_{\kappa+\widetilde{\tau}(({\Theta_0^2}/{2})n)}
- 2M_{\tau(({\Theta_0^2}/{2})n)} \bigr)^2\bigr]\\
&&\qquad\leq
\frac{5}{2}\mathbb{E}\biggl[\kappa+\widetilde{\tau}\biggl(\frac{\Theta
_0^2}{2}n\biggr) - \tau\biggl(\frac{\Theta_0^2}{2}n\biggr)\biggr]\\
&&\qquad= \frac{5}{3}\bigl((\mathbb{E}[\kappa] + \log n) - (\log n) \bigr)
+O(1)\\
&&\qquad= O(1) .
\end{eqnarray*}

Since
\[
\tau\biggl(\frac{\Theta^2_0}{2}n\biggr) \leq\sigma(n) \leq
\kappa+ \widetilde{\tau}\biggl(\frac{\Theta^2_0}{2}n\biggr)+\rho,
\]
and $\widetilde{\tau}(\frac{\Theta_0^2}{2}n)$ and $\tau(\frac
{\Theta_0^2}{2}n)$ differ only by
a quantity which has $O(1)$ second moment,
\begin{eqnarray*}
\mathbb{E}[\sigma(n)] &=& \mathbb{E}\biggl[\tau\biggl(\frac{\Theta
_0^2}{2}n\biggr)\biggr] + O(1) ,\\
\operatorname{Var}[\sigma(n)] &=& \operatorname{Var}\biggl[\tau\biggl
(\frac
{\Theta_0^2}{2}n\biggr)\biggr]
+ O\Biggl(\sqrt{\operatorname{Var}\biggl[\tau\biggl(\frac{\Theta
_0^2}{2}n\biggr)\biggr]}\Biggr) +
O(1) .
\end{eqnarray*}
Consequently, the theorem is proved as a consequence of the upper bound
asymptotics~(\ref{eq:upper-bound-asymptotics}) established at the
beginning of this proof.
\end{pf}

\subsubsection*{A Brownian digression}

Because we can compute $\operatorname{Var}[M_t]=\frac{5}{8}t$, we
can use
martingale central limit theorem ideas [\citet{Rebolledo-1980},
\citet{Whitt-2007}]
to show that
\[
\tau\approx\frac{2}{3}
\Biggl(
\log X_\tau+ \sqrt{7}B_\tau+2C - \log2
-\log\int^\infty_0\exp\Biggl(-\frac{3}{2} u +
\sqrt{\frac{5}{2}}\widetilde{B}_u^\tau\Biggr)\,\dd u
\Biggr)
,
\]
in the sense of weak convergence,
for $\widetilde{B}$ a standard Brownian motion not independent of $B$,
and $\sigma_\tau$ a stopping time for $B$ with expectation $\mathbb
{E}[\sigma_\tau]=\tau$. The distribution of the Dufresne integral
\[
\int^\infty_0\exp\Biggl(-\frac{3}{2} u + \sqrt{\frac
{5}{2}}\widetilde{B}_u^\tau\Biggr)\,\dd u
\]
is known explicitly [\citet{Dufresne-1990}, \citet
{Yor-1992a}]; however, its
contribution to the above is dominated by other terms.

\subsubsection*{Recovery of logarithmic excess result}%

Inspection of the growth process analysis shows that the logarithmic
excess occurs
before (say) time $n/\log n$, whereas our discussion of the maximum
lateral deviation
of $\partial\mathcal{C}(p^-,p^+)$ [with $\dist(p^-,p^+)=n$] shows
that the intersection of the two growth processes
occurs at an $x$-coordinate uniformly distributed over the range from
$p^-$ to~$p^+$.

This indicates that the asymptotic excess of the upper or the lower
semipe\-rimeter route for $\partial\mathcal{C}(p^-,p^+)$ should have
leading term $\frac{4}{3}\log n$, which agrees with the rigorous
arguments for the asymptotic behavior of the mean excess obtained in
\citeauthor{AldousKendall-2007} [(\citeyear{AldousKendall-2007}),
Theorem 3] and similar
higher-dimensional results obtained by \citeauthor{BoroczkySchneider-2008}
[(\citeyear{BoroczkySchneider-2008}), Theo-\break rem~1.3]
[compare the planar arguments of \citet{RenyiSulanke-1968}, Satz~5].

Note that the methods of the proof of Theorem~\ref{thm:variance-bound}
bear a family resemblance to the Markov chain methods which \citet
{Groeneboom-1988}
and \citet{CaboGroeneboom-1994} applied to problems concerning
convex hulls of Poisson point patterns.

\subsection{True network geodesics}\label{sec:true-geodesics}
We can now deduce that the two semipe\-rimeter routes provided by
$\partial\mathcal{C}(p^-,p^+)$ will often \textit{not} be network
geodesics. The boundary $\partial\mathcal{C}(p^-,p^+)$ is composed of
the initial parts of four independent growth processes, contributing
four independent initial excesses, each of mean $\frac{2}{3}\log n$\vadjust{\goodbreak}
and variance proportional to $\log n$; the remainder of the excess will
be of order less than $\sqrt{\log n}$. Accordingly, there is an even
chance that the least excess is achieved by crossing over from top side
to bottom side so as to use the smallest possible $\frac{2}{3}\log n
\pm\mbox{const.}\times\sqrt{\log n}$ contribution. Calculations of
the caricature of Section~\ref{sec:caricature} make it plain that such
a crossover can be achieved at the very modest price of adding just a
bounded term to the excess, and therefore there is a substantial
positive chance that one of the crossover routes is shorter than either
of the semiperimeter routes.

In fact, we conjecture that the two semiperimeter routes provided by
$\partial\mathcal{C}(p^-$, $p^+)$ are never network geodesics; in
particular, it should be possible to achieve modest reductions in the
excess by using crossovers very close to source and destination nodes
$p^-$ and $p^+$.


\subsubsection*{Lower bound for the Poissonian city}

Nevertheless, the semiperimeter routes supplied by $\partial\mathcal
{C}(p^-,p^+)$ are good approximations to true network geodesics; we
show this by establishing that their mean excess can be compared with
a lower bound on possible path lengths. Indeed, because we are working
in the specific situation of a Poisson line process, we can derive a
stronger and simpler version of the $\Omega(\sqrt{\log n})$ lower
bound argument of \citet{AldousKendall-2007}, Theorem 2.
\begin{theorem}\label{thm:poisson-lower-bound}
In the Poissonian city network, consider \textit{any} path from $p^-$ to
$p^+$ (in the sense described in Section~\ref{sec:connections}). If
$\dist(p^-,p^+)=n$, then the path must have mean excess exceeding
\[
2 \bigl(\log4 - \tfrac{5}{4}\bigr)\log n + o(\log n)
=
0.27258872\ldots\log n + o(\log n)
.
\]
\end{theorem}
\begin{pf}
Let $\mathcal{C}(\origin,+)$ be the cell containing the positive
$x$-axis of the tessellation formed from the Poisson line process $\Pi
$ by deleting all lines intercepting the positive $x$-axis. Consider
the vertical line $\ell_x$ through $(x,0)$ and let $-L_x^-$, $L_x^+$
be the distances along this line to $\partial\mathcal{C}(\origin,+)$
running down and up, respectively. Any network geodesic $\gamma$ from
$\origin$ to any other point $p$ on the positive $x$-axis, constructed
according to the recipe in Section~\ref{sec:connections}, must lie
between locations $-L^-$ and $L^+$ on $\ell_x$.
This is a consequence of the convexity of $\mathcal{C}(0,+)$.
Consider such a network geodesic, or indeed a general regular path
$\gamma$ lying within these bounds, and let $\theta_x\in(-\frac
{\pi}{2},\frac{\pi}{2})$ be the angle made with the horizontal by
$\gamma$ when encountering $\ell_x$ for the first time. If $\dist
(\origin,p)=n$, then the mean excess of $\gamma$ must exceed
\[
\mathbb{E}\biggl[\int_1^n(\sec\theta_x-1)\,\dd x\biggr]
\geq
\frac{1}{2}\int_1^n \mathbb{E}[\theta_x^2] \,\dd x
=
\int_1^n\int_0^{\pi/2}\mathbb{P}[|\theta_x|>u] u\,\dd u \,\dd x .
\]
Consider the probability of there being no lines of $\Pi$ which both
(a) hit $\ell_x$ in the range between $-L^-$ and $L^+$ signed
distances from the $x$-axis and (b) form an angle to the horizontal
which is less than $u$ in absolute value. The density of the angle to
the horizontal is $\frac{1}{2}\cos\theta$ for $-\frac{\pi
}{2}<\theta<\frac{\pi}{2}$, while the patterns of lines hitting
$\ell_x$ above and below the $x$-axis are independent. Consequently,
\[
\mathbb{P}[|\theta_x|>u]
\geq
\mathbb{E}\bigl[\exp\bigl(-(L^-_x+L^+_x)\sin u\bigr)\bigr]
\geq
(\mathbb{E}[\exp(- u L^+_x)])^2\vadjust{\goodbreak}
.
\]
Considerations from stochastic geometry determine the distribution of
$L^+_x$ and thus show that
\begin{eqnarray*}
\mathbb{E}[\exp(- u L^+_x)]
&=&
\int_0^1\mathbb{P}[\exp(- u L^+_x)>z]\,\dd z\\
&=&
\int_0^\infty\mathbb{P}[\exp(- u L^+_x)>e^{-s}]e^{-s}\,\dd s\\
&=&
1 - \int_0^\infty e^{-s}\mathbb{P}\biggl[L^+_x\geq\frac{s}{u}\biggr]\,\dd s\\
&=&
1 - \int_0^\infty\exp\Biggl(-s - \frac{1}{2}\Biggl(\sqrt
{x^2+\frac{s^2}{u^2}}-x\Biggr)\Biggr)\,\dd s
.
\end{eqnarray*}
Consequently,
\begin{eqnarray*}
&&\mathbb{E}\biggl[\int_1^n(\sec\theta_x-1)\,\dd x\biggr]\\
&&\qquad\geq
\int_1^n\int_0^{\pi/2}\mathbb{P}[|\theta_x|>u] u\,\dd u \,\dd x
\\
&&\qquad\geq
\int_1^n\int_0^{\pi/2}
\Biggl(1 - \int_0^\infty\exp\Biggl(-s - \frac{1}{2}\Biggl(\sqrt
{x^2+\frac{s^2}{u^2}}-x\Biggr)\Biggr)\,\dd s\Biggr)^2
u\,\dd u \,\dd x .
\end{eqnarray*}

Suppose that in the above we could control the error arising from
replacing $\sqrt{x^2+\frac{s^2}{u^2}}-x$ by its upper bound $\frac
{1}{2}s^2/(x u^2)$.
We would then need to estimate
\[
\int_1^n\int_0^{\pi/2}
\biggl(1 - \int_0^\infty\exp\biggl(-s - \frac{s^2}{4 u^2 x}
\biggr)\,\dd s\biggr)^2
u\,\dd u \,\dd x .
\]
However, we can in fact estimate
\[
\int_0^\infty\exp\biggl(-s - \frac{s^2}{2p^2}\biggr)\,\dd s
=
p\exp\biggl(\frac{p^2}{2}\biggr)\int_p^\infty\exp\biggl(-\frac
{s^2}{2}\biggr)\,\dd s
\]
using classical results on Mill's ratio, namely the excellent upper
bound of \citet{Sampford-1953} [see also \citet
{Baricz-2008} for a treatment based on a monotone form of l'H\^opital's rule]:
%
%
\begin{equation}\label{eq:sampford}
\exp\biggl(\frac{p^2}{2}\biggr)\int_p^\infty\exp\biggl(-\frac
{s^2}{2}\biggr)\,\dd s
\leq
\frac{4}{\sqrt{p^2+8}+3p}\qquad
\mbox{for }p>-1 .
\end{equation}
Accordingly, setting $v=\sqrt{x}u$ and letting $n$ tend to $\infty$,
we obtain
\begin{eqnarray*}
&&\int_1^n\int_0^{\pi/2}
\biggl(1 - \int_0^\infty\exp\biggl(-s - \frac{s^2}{4 u^2 x}
\biggr)\,\dd s\biggr)^2
u\,\dd u \,\dd x
\\[-3pt]
&&\qquad=\int_1^n\int_0^{\sqrt{x}\pi/2}
\biggl(1 - \int_0^\infty\exp\biggl(-s - \frac{s^2}{4 v^2}
\biggr)\,\dd s\biggr)^2
v\,\dd v\,\frac{\dd x}{x}
\\[-3pt]
&&\qquad\geq
\int_1^n\int_0^{\sqrt{x}\pi/2}
\biggl(
\frac{\sqrt{v^2+4}-v}{\sqrt{v^2+4}+3v}
\biggr)^2
v\,\dd v\,\frac{\dd x}{x}
\\[-3pt]
&&\qquad\sim
\int_1^n\int_0^\infty
\biggl(
\frac{\sqrt{v^2+4}-v}{\sqrt{v^2+4}+3v}
\biggr)^2
v\,\dd v \,\frac{\dd x}{x}\\[-3pt]
&&\qquad=
\biggl(\log4 - \frac{5}{4}\biggr)\log n.
\end{eqnarray*}
Here, the $v$-integral is evaluated by making the substitutions
$v=\sinh t$ and then $y=2-e^{-2t}$.

The proof is completed by bounding the error arising from the approximation
of $\sqrt{x^2+s^2/u^2}-x$ by $\frac{1}{2}s^2/(x u^2)$:
\begin{eqnarray*}
&&\int_1^n\int_0^{\pi/2} \Biggl(1 - \int_0^\infty\exp\Biggl(-s -
\frac{1}{2}\Biggl(\sqrt {x^2+\frac{s^2}{u^2}}-x\Biggr)\Biggr)\,\dd
s\Biggr)^2
u\,\dd u \,\dd x\\[-3pt]
&&\qquad= \int_1^n\int_0^{\pi/2} \Biggl( 1 - \int_0^\infty
e^{-s}\Biggl( \exp\Biggl(-
\frac{1}{2}\Biggl(\sqrt{x^2+\frac{s^2}{u^2}}-x
\Biggr)\Biggr)\\[-3pt]
&&\qquad\quad\hspace*{157.1pt}{} - \exp\biggl(- \frac{s^2}{4 u^2
x}\biggr) \Biggr)
\,\dd s\\[-3pt]
&&\qquad\quad\hspace*{126.8pt}{} - \int_0^\infty \exp\biggl(-s -
\frac{s^2}{4 u^2 x}\biggr) \,\dd s \Biggr)^2
u\,\dd u \,\dd x \\[-3pt]
&&\qquad\geq \int_1^n\int_0^{\pi/2} \biggl(1 -
\int_0^\infty\exp\biggl(-s - \frac{s^2}{4 u^2 x} \biggr)\,\dd
s\biggr)^2
u\,\dd u \,\dd x\\[-3pt]
&&\qquad\quad\hspace*{0pt}{}- 2 \int_1^n\int_0^{\pi/2} \biggl(1 -
\int_0^\infty\exp\biggl(-\overline{s} - \frac {\overline{s}^2}{4 u^2
x}\biggr)\,\dd\overline{s}\biggr)
\\[-3pt]
&&\qquad\quad\hspace*{58.6pt}{}\times \int_0^\infty e^{-s} \Biggl(
\exp\Biggl(-\frac{1}{2}\Biggl(\sqrt{x^2+\frac{s^2}{u^2}}-x
\Biggr)\Biggr)\\[-3pt]
&&\qquad\quad\hspace*{163.1pt}{}- \exp\biggl(-\frac{s^2}{4 u^2
x}\biggr) \Biggr)\,\dd s\, u\,\dd u \,\dd x .
\end{eqnarray*}
We need to bound the second term, and we do this by invoking
Birnbaum's (\citeyear{Birnbaum-1942}) very good lower bound\vadjust{\goodbreak}
on Mill's ratio. After some manipulation this yields
\[
1 - \int_0^\infty\exp\biggl(-\overline{s} - \frac{\overline
{s}^2}{4 u^2 x}\biggr)\,\dd\overline{s}
\leq
\frac{\sqrt{u^2 x+2}-u\sqrt{x}}{\sqrt{u^2 x+2}+u\sqrt{x}} .
\]
It thus suffices to bound
\begin{eqnarray*}
&&2\int_1^n\int_0^{\pi/2} \frac{\sqrt{u^2 x+2}-u\sqrt{x}}{\sqrt{u^2
x+2}+u\sqrt{x}} \\
&&\qquad\hspace*{24.6pt}{}\times\int_0^\infty e^{-s} \bigl(
e^{-({1}/{2})(\sqrt{x^2+{s^2}/{u^2}}-x)} - e^{-{s^2}/({4 u^2
x})} \bigr)
\,\dd s \,u\,\dd u \,\dd x\\
&&\qquad\leq \int_1^n\int_0^{\pi/2} \frac{2}{\sqrt{{u^2
x}+2}+u\sqrt{x}} \\
&&\qquad\hspace*{50pt}{}\times\int_0^\infty e^{-s} \bigl(
e^{-({1}/{2})(\sqrt{x^2+{s^2}/{u^2}}-x)} - e^{-{s^2}/({4 u^2
x})} \bigr)
\,\dd s \,\dd u \,\frac{\dd x}{\sqrt{x}}\\
&&\qquad\leq \sqrt{2} \int_1^\infty\int_0^{\pi/2}\int_0^\infty e^{-s}
\min\Biggl\{1, \frac{s^2}{4u^2 x} - \frac{1}{2}\Biggl(\sqrt{x^2+\frac{s^2}{u^2}}-x\Biggr)
\Biggr\}
\,\dd s \,\dd u \,\frac{\dd x}{\sqrt{x}}\\
&&\qquad\leq \sqrt{2} \int_0^\infty\int_0^{\pi/2}\int_0^\infty e^{-s}
\min\biggl\{1, \frac{1}{16}\frac{s^4}{u^4 x^3} \biggr\} \,\dd s \,\dd u \,\frac{\dd
x}{\sqrt{x}} ,
\end{eqnarray*}
where we use the fact that over the range $0\leq p<\infty$, the
function $p\mapsto e^{-p}$ is nonnegative, decreasing and Lipschitz
with Lipschitz constant $1$,
and also that $1+\frac{1}{2}p-\sqrt{1+p}\leq\frac{1}{8}p^2$ if
$0\leq p<\infty$ (use finite Taylor series expansion).

We split the $x$-integral at $x^3=s^4/(16u^4)$; the first part is
bounded by
\begin{eqnarray*}
&&\sqrt{2} \int_0^\infty\int_0^{\pi/2} \int_0^{(s^4/(16u^4))^{1/3}}
e^{-s} \min\biggl\{1, \frac{1}{16}\frac{s^4}{u^4 x^3} \biggr\}\, \frac{\dd
x}{\sqrt{x}} \,\dd u \,\dd s
\\
&&\qquad= \sqrt{2} \int_0^\infty\int_0^{\pi/2}
\int_0^{(s^4/(16u^4))^{1/3}} e^{-s}\, \frac{\dd x}{\sqrt{x}}
\,\dd u \,\dd s\\
&&\qquad= 3\sqrt{2} \pi^{1/3}\Gamma\biggl(\frac{5}{3}\biggr) ,
\end{eqnarray*}
while the second part is bounded by
\begin{eqnarray*}
&&\sqrt{2} \int_0^\infty\int_0^{\pi/2}
\int_{(s^4/(16u^4))^{1/3}}^\infty e^{-s} \min\biggl\{1,
\frac{1}{16}\frac{s^4}{u^4 x^3} \biggr\} \,\frac{\dd x}{\sqrt{x}} \,\dd u \,\dd
s\\
&&\qquad= \sqrt{2} \int_0^\infty\int_0^{\pi/2}
\int_{(s^4/(16u^4))^{1/3}}^\infty e^{-s} \frac{1}{16}\frac{s^4}{u^4
x^3} \,\frac{\dd x}{\sqrt{x}} \,\dd u \,\dd s
\\
&&\qquad = \frac{3}{5} \sqrt{2}\pi^{1/3}\Gamma\biggl(\frac{5}{3}\biggr) .
\end{eqnarray*}
This establishes the desired result since we may apply the lower bound
on mean excess path length to (a) the half of the geodesic running from
$p^-$ to midway and (b) the other half running from midway to $p^+$.
\end{pf}

The results of this section justify the focus in the remainder of this
paper on the semiperimeter routes provided by $\partial\mathcal{C}(p^-,p^+)$:
while semiperimeter routes do differ from network geodesics,
their use nevertheless does not incur a great penalty; they are
produced by a geometric algorithm which is certainly unsophisticated,
but, on the other hand, is explicit, and they are amenable to exact
calculations.

\section{Traffic flow in the Poissonian city}\label{sec:flow}
We now consider traffic flow in the network produced by this Poisson
line process. To do this, we first compute the mean flow through a line
at the center of the disc. More precisely, we condition on there being
a (horizontal) line of the line process running through the origin
$\origin$ and consider the flow through $\origin$ which results if
every pair of $x$ and $y$ in $\ball(\origin,n)$ generates an
infinitesimal flow of amount ${d}x \,{d}y$ divided equally between the
%
%
\begin{figure}

\includegraphics{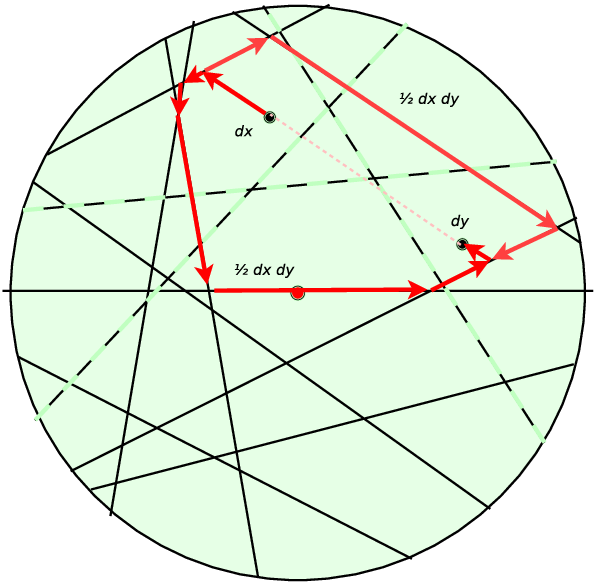}

\caption{Illustration of the flow generated between two points in a
``Poissonian
city.''}\label{fig:flow}
\end{figure}
two possible routes given by the semiperimeter algorithm (see Figure
\ref{fig:flow}).

\subsection{First order calculations at the center}\label{sec:first-order}

Recall~(\ref{eq:region}) from the \hyperref
[sec:introduction]{Introduction}: the flow
through the center is measured by the $4$-volume of $\mathcal{D}_n$, where
\[
\mathcal{D}_n =
\{
(p^-, p^+)\in\ball(\origin,n)^2 \dvtx p^-_1<p^+_1, \origin\in
\partial\mathcal{C}(p^-,p^+)
\} ,
\]
and we seek to understand the large-$n$ statistical behavior of this $4$-volume.
Indeed (bearing in mind that we have conditioned on there being a line
through $\origin$),
the distribution of the total traffic through the origin $\origin$ is
given by the distribution of
%
%
\begin{eqnarray}\label{eq:total-flow}
T_n &=& 
\frac{1}{2}\iint\mathbb{I}_{[(p^-,p^+)\in\mathcal{D}_n]}\,\dd
p^-\,\dd p^+
\nonumber\\[-8pt]\\[-8pt]
&=& 
\frac{1}{2}\iint_{\ball(\origin,n)^2}\mathbb{I}_{[p^-_1<p^+_1,
\origin\in\partial\mathcal{C}(p^-,p^+)]}\,\dd p^-\,\dd p^+ .
\nonumber
\end{eqnarray}


The mean can be obtained asymptotically using direct arguments.
\begin{theorem}\label{thm:first-moment}
The mean flow through a line at the center is given by
%
%
\begin{equation}\label{eq:mult-int-1}
\mathbb{E}[T_n] =
\int_0^\pi\int_0^n\int_0^n \exp\biggl(-\frac{1}{2}(r+s-\rho
)\biggr) r\,\dd r \,s\,\dd s \,\theta\,\dd\theta,
\end{equation}
where $\rho=\sqrt{r^2+s^2+2 r s \cos\theta}$. Asymptotically, as
$n\to\infty$,
%
%
\begin{equation}\label{eq:first-moment}
\mathbb{E}[T_n] \sim2 n^3 .
\end{equation}
\end{theorem}
\begin{pf}
Equation~(\ref{eq:mult-int-1}) follows from simple stochastic geometry
of Poisson line processes (illustrated in Figure~\ref{fig:expect}),\vadjust{\goodbreak}
taking care not to double-count flow between the unordered points $p^+$
and $p^-$. Note that when $p^+$ and $p^-$ are on opposing sides of the
line conditioned to hit the origin, none of the flows between these two
points will run through the origin. Indeed, mean flows between points
$p^+$, $p^-$ in the upper half-plane will account for exactly half the
total mean flow through the origin.

%
%
\begin{figure}

\includegraphics{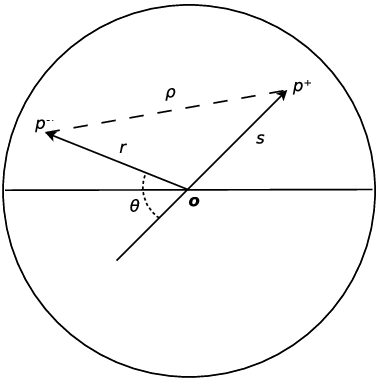}

\caption{Illustration of the geometry
represented by the multiple integral in (\protect\ref{eq:mult-int-1}) using
$r$, $s$, $\theta$. The segment $p^-p^+$ is not separated from
the origin $\origin$ exactly when no lines of the line process pass
through both of $\origin p^-$ and
$\origin p^+$.}\label{fig:expect}
\end{figure}

The asymptotics follow by application of scaling by $n$, the symmetry
between $s$ and $r$, and the inequality $\sqrt{1-2z}\leq1-z$ for
$z\geq\frac{1}{2}$. Indeed,
\begin{eqnarray*}
\mathbb{E}[T_n] &=& 2n^4 \int_0^\pi\int_0^1\int_0^s \exp
\biggl(-\frac{n}{2}(r+s-\rho)\biggr) r\,\dd r \,s\,\dd s \,\theta\,\dd
\theta\\
&=& 2n^4 \int_0^\pi\int_0^1\int_0^1 \exp\biggl(-\frac{n
s}{2}\bigl(r+1-\sqrt{r^2+1+2 r \cos\theta}\bigr)\biggr) r\,\dd r \,s^3\,\dd s\,
\theta\,\dd\theta.
\end{eqnarray*}

The region of the integral corresponding to
$\int_{\theta_n}^{\pi/2}\int_0^1\int_0^1$ (for $\theta_n>0$) is
bounded above by
\begin{eqnarray*}
&&2n^4 \int_{\theta_n}^\pi\int_0^1\int_0^1 \exp\biggl(-\frac{n
s}{2}\bigl(r+1-\sqrt{r^2+1+2 r \cos\theta}\bigr)\biggr) r\,\dd r\, s^3\,\dd s
\theta\,\dd\theta\\
&&\qquad\leq
2n^4 \int_{\theta_n}^\pi\int_0^1\int_0^1 \exp\biggl(-\frac{n s
r}{2(r+1)}(1-\cos\theta)\biggr) r\,\dd r \,s^3\,\dd s \,\theta\,\dd\theta\\
&&\qquad\leq
(\pi-\theta_n)^2 n^4 \int_0^1\int_0^1 \exp\biggl(-\frac{n s
r}{2(r+1)}(1-\cos\theta_n)\biggr) r\,\dd r \,s^3\,\dd s\\
&&\qquad\leq
(\pi-\theta_n)^2 n^4 \int_0^1\int_0^\infty\exp\biggl(-\frac{n s
r}{4}(1-\cos\theta_n)\biggr) r\,\dd r \,s^3\,\dd s\\
&&\qquad=
(\pi-\theta_n)^2 \biggl(\frac{4}{1-\cos\theta_n}\biggr)^2 n^2
\int_0^1 s\,\dd s=
8 \biggl(\frac{\pi-\theta_n}{1-\cos\theta_n}\biggr)^2 n^2.
\end{eqnarray*}

Consider the region $\int_0^{\theta_n}\int_0^1\int_0^1$. Using a
Taylor series expansion of $1-\sqrt{1-2z}$ and the approximation
$\theta/\sin\theta\searrow1$ as $\theta\searrow0$ (so long as
$0<\theta<\pi/2$), we deduce that
\begin{eqnarray*}
\mathbb{E}[T_n] &\sim&
2 n^4 \int_0^1\int_0^1 \int_0^\infty\exp\biggl(-\frac
{nsr}{2(1+r)}u\biggr)\,\dd u \,s^3\,\dd s \,r\,\dd r\\
&=&
2 n^4 \int_0^1\int_0^1 \frac{2(1+r)}{n s r} s^3\,\dd s\, r\,\dd r
= 2 n^3 .
\end{eqnarray*}
\upqed\end{pf}

Taking some extra care over the analysis, it is possible to bound the
order of the error in~(\ref{eq:first-moment}). We state this
without proof.
\begin{cor}\label{cor:sharper-asymptotic}
The asymptotic in Theorem~\ref{thm:first-moment} can be sharpened to
\[
\mathbb{E}[T_n] \sim2 n^3 + O\bigl(n^2\sqrt{n}\bigr) \qquad\mbox{as
}n\to\infty.
\]
\end{cor}
%

\subsection{Mean flow averaged over entire disc}\label
{sec:location-dependence}

Regional variation of expected flow over the disc is to be expected:
flow at the boundaries should be lower than at the center. Indeed, one
can calculate the mean flow per unit length in the network as follows.

The mean total length of the intersection of the Poisson line pattern
with the disc is given by (mean number of lines hitting disk) $\times$
average intersection length:
\[
(2\pi n)
\times\biggl(
\frac{1}{2n}\int_{-n}^n 2\sqrt{n^2-x^2}\,\dd x\biggr)
=
\frac{\pi^2 n^2}{2} ,
\]
and this is therefore the mean total network length in the disc.

On the other hand, the mean Euclidean distance between two independent
uniformly random points in the disc is given by
\begin{eqnarray*}
&&\frac{1}{\pi n^2}\int_0^{2\pi}\int_0^n
\frac{1}{\pi n^2}\int_0^{2\pi}\int_0^n
\sqrt{u^2+v^2-2 u v \cos(\alpha-\beta)}
v\,\dd v \,\dd\beta\,
u\,\dd u \,\dd\alpha
\\
&&\qquad=
\frac{8 n}{5\pi}\int_0^\pi\int_0^1
\sqrt{u^2+1-2 u \cos\theta}
u\,\dd u \,\dd\theta
\\
&&\qquad=
\frac{8 n}{5\pi}\int_0^{\pi/2}\int_0^{2\cos\phi}
s^2
\,\dd s \,\dd\phi
=
\frac{128}{45\pi} n ,
\end{eqnarray*}
where the first step uses various symmetries and rescaling, and the
second step changes to polar coordinates based at $u=1$ and $\theta=0$
in an implicit use of Crofton's method. [This calculation is a special
case of a classic calculation in geometric probability, surveyed in
\citet{Santalo-1976}, Chapter 12.7 Note (6).] By the previous results on
lengths of network geodesics, the mean \textit{network} distance differs
only by an extra logarithmic contribution.

Hence, the mean flow per unit length, if each pair of points exchanges
just one infinitesimal unit of traffic and this is averaged over the
network, is asymptotic to
\[
\frac{1}{2} \frac{(\pi n^2)^2 \times
({128}/({45\pi})) n}{{\pi^2 n^2}/{2}}
=
\frac{128}{45\pi} n^3
=
1.9052\ldots n^3 .
\]
This analysis does not take account of routes which move outside the
perimeter of the disc; however, the effect of these routes can be shown
to be negligible. [The key observation is based on Theorem \ref
{thm:lateral-displacement}: if both source and destination nodes $p^\pm
$ are at least $2\sqrt{(1+\varepsilon)n\log n}$ from the perimeter of
a disc of radius $n$, and $n$ is large enough, then points outside the
disc have probability at most $O(n^{-(1+\varepsilon)})$ of lying
within $\mathcal{C}(p^-,p^+)$. Thus, mean total length outside the
disc is a boundary rather than an area effect.]
In conclusion, and unsurprisingly, mean flow over a typical line is
slightly smaller than mean flow over a line at the center of the disc.

\subsection{An improper anisotropic limiting line process}\label
{sec:limiting-object}

We can represent the scaling limit of the distribution of traffic flow
through the center of the Poissonian city by using an improper
stationary anisotropic Poisson line process.

\begin{figure}

\includegraphics{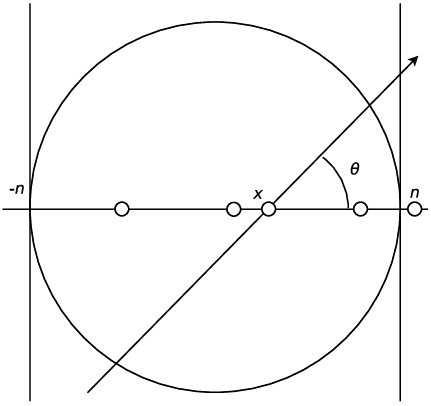}

\caption{The Poisson line process can be
represented in terms of a Poisson process of points scattered along the
$x$-axis, through each of which there runs a line making an angle
$\theta\in(0,\pi)$ with the $x$-axis, with density
$\frac{1}{2}\sin\theta$. Here, we show this against a backdrop of the
disc $\ball(\origin,n)$.}\label{fig:limit1}
\end{figure}

We use the alternate coordinatization of a unit rate stationary
isotropic Poisson line process, as in Figure~\ref{fig:limit1}, using
coordinate $x$ for the intersection along the $x$-axis and $\theta$
for line direction. The $x$-axis intersections then form a stationary
Poisson point process, while the angle density is $\frac{1}{2}\sin
\theta$ for $0<\theta<\pi$.

%
%

%
\begin{figure}

\includegraphics{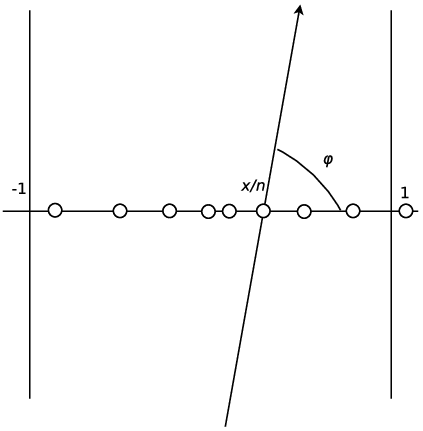}

\caption{This illustrates the result of scaling by $1/n$
in the $x$-direction and $1/\sqrt{n}$ in the
$y$-direction.}\label{fig:limit2}\vspace*{-3pt}
\end{figure}

Rescale to shrink the $x$-axis by a factor of $1/n$, so $\widetilde
{x}=x/n$. Guided by previous results, shrink the $y$-axis by a
different amount, $1/\sqrt{n}$, so $\widetilde{y}=y/\sqrt{n}$. Thus,
$\theta$ is transformed into a new angle $\phi$ (see Figure \ref
{fig:limit2}), where
\[
\tan\phi= \sqrt{n}\tan\theta.
\]
In the new coordinates of $\widetilde{x}$ and $\phi$, the line
process can be parametrized as a nonstationary Poisson point process on
$\widetilde{x}\dvtx\phi$ space with intensity
\[
\frac{1}{2} \frac{\tan\phi\sec^2\phi}{(1+({1}/{n})\tan
^2\phi)^{3/2}}\,\dd\widetilde{x} \,\dd\phi
\nearrow
\frac{1}{2}\tan\phi\sec^2\phi\,\dd\widetilde{x} \,\dd\phi.
\]
We can represent this as a coupling construction: based on an improper
stationary anisotropic Poisson line process with intensity $\frac
{1}{2}\tan\phi\sec^2\phi\,\dd\widetilde{x} \,\dd\phi$ in
$\widetilde
{x}\dvtx\phi$ coordinates, we can achieve a proper stationary \textit
{isotropic} Poisson line process at scale $n$ by randomly thinning the
lines with retention probability depending monotonically on the line slope.

\begin{figure}

\includegraphics{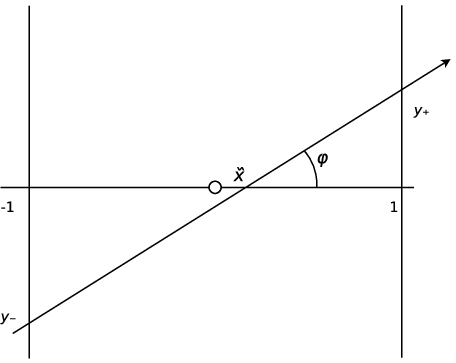}

\caption{This illustrates the parametrization of lines from
the improper limiting line process as intercepts on two parallel
$y$-axes.}\label{fig:limit3}\vspace*{-3pt}
\end{figure}

Moreover, this limiting object may be cleanly represented using a
further set of coordinates. Represent each line of the line process by
its intercepts $y_+$ and $y_-$ on the vertical axes $x=1$ and $x=-1$
(see Figure~\ref{fig:limit3}). The intensity then becomes
\[
\tfrac{1}{4}\,\dd y_+\,\dd y_- .
\]
In particular, while the new improper Poisson line process is
anisotropic, it nevertheless does possess special affine shear
symmetries, namely the symmetries produced by those area-preserving
linear transformations which leave all vertical axes invariant.

%


This construction enables us to identify the limiting behavior for
$T_n$, as follows.\vspace*{-3pt}
\begin{theorem}\label{thm:scaling-limit}
The scaled quantity
$T_n/n^3$
has a limiting distribution given by the analogous flow at the center
for the limiting improper stationary anisotropic Poisson line process
given above.\vspace*{-3pt}
\end{theorem}

Indeed, we can relate scaled finite-$n$ instances to the limiting case
by a coupling argument involving the addition of further lines;
however, the resulting\vadjust{\goodbreak} almost sure limit is not monotonically
decreasing since it will involve a double integral [as in~(\ref
{eq:total-flow})] taken over ever-increasing regions of the vertical
strip in Figure~\ref{fig:limit3}.

Before proving this theorem, we show that the mean flow at the center
for this limit is in agreement with the asymptotics given in Theorem
\ref{thm:first-moment}.\vspace*{-3pt}
\begin{lem}\label{thm:scaling-limit-mean}
The flow at the center for the limiting improper stationary anisotropic
Poisson line process given above has mean value $2$.\vspace*{-3pt}
\end{lem}
\begin{pf}
Consider first the probability that the line segment from $(-a,t)$ to
$(u,t)$ is not separated from the origin by the improper line process.
The mean measure of lines implementing such a separation (measured
using the intensity measure of the improper process) is $A+B+C$, where
the contribution $A$ arises from separating lines hitting the
upper-left shaded triangle in Figure~\ref{fig:improper-mean}, the
contribution $C$ arises from those hitting the upper-right shaded
triangle in Figure~\ref{fig:improper-mean} and $B$ is derived from the
contribution of the remaining separating lines.

%
%
\begin{figure}

\includegraphics{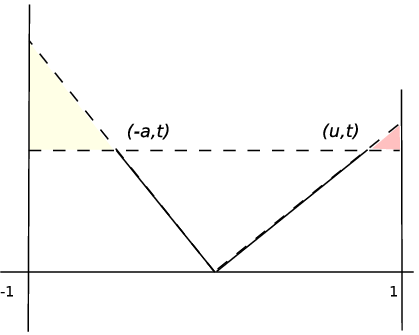}

\caption{Computing the mean flow through the center for the
flow based on the improper anisotropic line process
limit.}\label{fig:improper-mean}
\end{figure}

Then,
\[
A =
\frac{1}{4}\int_t^{t/a}
\biggl(\biggl(s - 2\frac{s-t}{1-a}\biggr)+s\biggr)
\,\dd s
=
\frac{t^2}{4}\biggl(\frac{1}{a}-1\biggr)
\]
and, similarly, $C=\frac{t^2}{4}(\frac{1}{u}-1)$. Finally,
\[
B = \frac{1}{4}\int_{-t}^t(t+s)\,\dd s = 2 \frac
{t^2}{4} .
\]
Consequently,
\[
A+B+C = \frac{t^2}{4}\biggl(\frac{1}{a}+\frac{1}{u}
\biggr) .\vadjust{\goodbreak}
\]
%
Since we are dealing with a Poisson process,
the required probability of the line segment from $(-a,t)$ to $(u,t)$
not being separated from the origin is
%
%
\begin{equation}\label{eq:affine-invariant}
\exp\biggl(-\frac{t^2}{4}\biggl(\frac{1}{a}+\frac{1}{u}
\biggr)\biggr) .
\end{equation}

Consider the special affine shear symmetries which leave the $x=0$ axis
fixed. Because of this symmetry group, it follows that the probability
of the line segment from $(-a,b)$ to $(u,v)$ not being separated from
the origin agrees with~(\ref{eq:affine-invariant}) when the line
segment from $(-a,b)$ to $(u,v)$ passes through $(0,t)$.
This occurs when $t=\frac{bu+av}{a+u}$; moreover, if we set $s=b-v$,
then $\dd t \,\dd s=\dd b \,\dd v$. Accordingly, the mean $4$-volume of the
region representing the flow through the center is given by
\begin{eqnarray*}
&&\int_0^1\int_0^1\int_0^\infty\int_{-(({a+u})/{a})t}^{
(({a+u})/{u})t}\exp\biggl(-\frac{t^2}{4}\biggl(\frac{1}{a}+\frac
{1}{u}\biggr)\biggr)\,\dd s \,\dd t \,\dd a \,\dd u
\\
&&\qquad=
\int_0^1\int_0^1\int_0^\infty(a+u)\biggl(\frac{1}{a}+\frac
{1}{u}\biggr)\exp\biggl(-\frac{t^2}{4}\biggl(\frac{1}{a}+\frac
{1}{u}\biggr)\biggr) t \,\dd t \,\dd a \,\dd u
\\
&&\qquad=
2 \int_0^1\int_0^1 (a+u) \,\dd a \,\dd u= 2 .
\end{eqnarray*}
\upqed\end{pf}
\begin{pf*}{Proof of Theorem~\ref{thm:scaling-limit}}
Consider the affine shear transformation $\mathcal{T}_n\dvtx[-1$, $1]\times
(0,\infty)\to[-1,1]\times(0,\infty)$ given by $\mathcal
{T}_n(u,v)=(n u,\sqrt{n} v)$. Define coupled random functions
\begin{eqnarray*}
&\displaystyle I_n\dvtx\bigl([-1,0]\times(0,\infty)\bigr)\times\bigl([0,1]\times
(0,\infty)\bigr) \to\{0,1\} ,&\\
&\displaystyle I_n(p,q) = \mathbb{I}_{[\mathcal{T}_n p\in\ball(\origin
,n)]}\mathbb{I}_{[\mathcal{T}_n q\in\ball(\origin,n)]}
\mathbb{I}_{[\origin\in\mathcal{C}(\mathcal{T}_n p,\mathcal{T}_n
q)]} .&
\end{eqnarray*}
So, $I_n$ depends implicitly on the underlying Poisson line process:
the previously described coupling construction shows that we can couple
different Poisson line processes for each $n$ so as to arrange that
$I_n(p,q)\to I(p,q)$ almost surely for Lebesgue almost all $p$, $q$,
where $I(p,q)$ is given by an analogous construction based on the
limiting improper anisotropic Poisson line process, and not using
$\mathcal{T}_n$.
Moreover, we can realize $T_n$ using
\[
T_n/n^3 = \frac{1}{2}\iint I_n(p,q) \,\dd p \,\dd q .
\]

From Theorem~\ref{thm:first-moment} and Lemma \ref
{thm:scaling-limit-mean}, it follows that
%
\[
\mathbb{E}\biggl[\frac{1}{2}\iint I_n(p,q) \,\dd p \,\dd q\biggr] \to
\mathbb{E}\biggl[\frac{1}{2}\iint I(p,q) \,\dd p \,\dd q\biggr] = 2 .
\]

On the other hand, if we restrict consideration to the finite measure
space $\Omega\times([-1,0]\times(0,K))\times
([0,1]\times(0,K))$ for any fixed $K$, then we may deduce
$L^1$-convergence of $I_n$ to $I$ via the dominated convergence theorem
since the indicator functions $I_n$ are bounded. [Here, $(\Omega
,\mathfrak{F},\mathbb{P})$ is the underlying probability space.]

It then follows from nonnegativity of the $I_n$, $I$ that we can apply
Fatou's lemma to deliver $L^1$-convergence on all of $\Omega\times
([-1,0]\times(0,\infty))\times([0,1]\times
(0,\infty))$
and so can deduce convergence in distribution as required:
\[
T_n/n^3 = \frac{1}{2}{\iint I_n(p,q) \,\dd p \,\dd q}
\mathop{\to} _{\mathcal{D}} \frac{1}{2}{\iint
I(p,q) \,\dd p \,\dd q,}
\]
viewed as random variables (functions of $\omega\in\Omega$).
\end{pf*}

Note that this proof also establishes uniform integrability of the
sequence of random variables $\{T_n/n^3\dvtx n\geq1\}$.

It is apparent from this construction that the limiting distribution is
largely insensitive to modest variations in the geometry of the city
[$\ball(\origin, n)$, or square of side $2n$, or$\ldots$];
however, we will not explore this here.


In principle it is possible that the limiting distribution of $T_n/n^3$
might be degenerate. That this is not the case follows rapidly from
representation of the limit in terms of the improper anisotropic
Poisson line process.
\begin{cor}\label{cor:non-degenerate}
The limiting distribution of $T_n/n^3$ is nondegenerate.
\end{cor}
\begin{pf}
Let $E_k$ be the event
\[
E_k = \biggl[\mbox{there is a line connecting }\{-1\}\times
\biggl[0,\frac{1}{k}\biggr]\mbox{ to }\{+1\}\times\biggl[0,\frac{1}{k}\biggr]\biggr] .
\]
Then, $E_k$ has positive probability for the improper anisotropic
Poisson line process; moreover, $E_1$, $E_2,\ldots$ form a
monotonically decreasing sequence of events whose intersection is a
null set. It follows from elementary measure theory that
\[
\mathbb{E}\biggl[\frac{1}{2} \iint I(p,q)\,\dd p \,\dd q ; E_k\biggr] \to
0 .
\]
However, simple constructions show positivity of the conditional expectation
\[
\mathbb{E}\biggl[\iint I(p,q)\,\dd p \,\dd q \Big| E_k\biggr]>0
\]
for each $k$. Because each event $E_k$ is of positive probability, it
follows that for each $\varepsilon>0$, we can find $k$ such that
\[
0 < \mathbb{E}\biggl[\frac{1}{2} \iint I(p,q)\,\dd p \,\dd q ; E_k\biggr] <
\varepsilon.
\]
It follows that the random variable
\[
\frac{1}{2} \iint I(p,q)\,\dd p \,\dd q
\]
cannot be deterministic,
and this establishes nondegeneracy of the limiting distribution of $T_n/n^3$.
\end{pf}


\subsection{Higher order moments}\label{sec:uniform-integrability}

\begin{figure}

\includegraphics{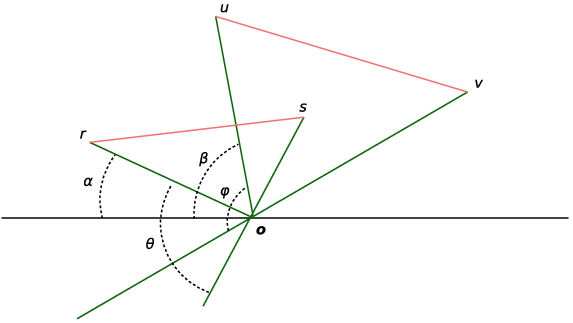}

\caption{Illustration of the event
$E(r,s,\theta,\alpha,u,v,\phi,\beta)$, which happens if neither of
the segments $u\to v$, $r\to s$ is separated from the origin by the
line process.}\label{fig:eight-fold}
\end{figure}

We have established that one can produce a coupling construction to
show that $T_n/n$ converges almost surely, and indeed in mean value, to
the corresponding quantity for the improper stationary anisotropic line
process. In fact, it is possible to establish convergence of moments of
order $2-\epsilon$ for $\epsilon\in(0,2)$; computations show that
the second moment $\mathbb{E}[T_n^2]$ is bounded by $\mathrm
{const.}\times n^6$, hence a uniform integrability argument may be
applied. The tiresomely complicated computations are omitted; we simply
indicate the general approach. The second moment can be expressed using
an eight-fold integral:
\begin{eqnarray*}
&&\int_0^n\int_0^n\int_0^\pi\int_0^\theta\int_0^n\int_0^n\int
_0^\pi\int_0^\phi
\mathbb{P}[E(r,s,\theta,\alpha,\\
&&\hspace*{158.3pt}u,v,\phi,\beta)]
\,\dd\beta\,\dd\phi\,
u\,\dd u \,v\,\dd v
\,\dd\alpha\,\dd\theta\,
r\,\dd r \,s\,\dd s ,
\end{eqnarray*}
where $E(r,s,\theta,\alpha,u,v,\phi,\beta)$ is the event that
neither of two line segments [$(r,\alpha)$-$(s,\alpha+\pi-\theta)$
and $(u, \beta,)$-$(v,\beta+\pi-\phi)$ when written in polar
coordinates] is separated from the origin by the line process (see
Figure~\ref{fig:eight-fold}).

%

Analysis of this eight-fold integral is complicated because the probability
\[
\mathbb{P}[E(r,s,\theta,\alpha,u,v,\phi,\beta)]
\]
takes on around eight different analytical forms, according to the
relative geometry of the line segments. Case-by-case mathematical
analysis shows that the multiple integral is bounded by $\mathrm
{const.}\times n^6$.

\section{Comparison with a Manhattan city}\label{sec:grid-structure}%
It is natural to ask how the Poissonian city might compare with an
alternate Manhattan city based on a more conventional Cartesian grid
structure of roads (a ``grid city''). Consider the case of a disc of
radius $n$ furnished with roads arranged in a fixed unit-length grid
structure. Suppose we wish to connect from the point at $-(u,v)$ to the
point at $(x,y)$ (using Cartesian coordinates). If $u$, $v$, $x$ and
$y$ are all nonnegative, then there is a wide variety of possible
network geodesic connections. Working with the most direct analogy to
the results described in Section~\ref{sec:flow}, suppose that traffic
from $-(u,v)$ to $(x,y)$ divides equally between the two extreme
network geodesics running from $-(u,v)$ to $(x,y)$. Working to $O(n^3)$
(allowing us to ignore some double-counting), the total traffic through
$\origin$ is made up of four contributions, arising from (i) $u=0,
v>0$, (ii) $u>0, v=0$, (iii) $x=0, y>0$, (iv) $x>0, y=0$. Each case
individually contributes a term of the form
\[
2 \mathop{\mathop{\sum\sum}_{(x,y)>0}}_{x^2+y^2\leq n^2}
\frac{1}{2} n
= \frac{\pi}{4} n^3 + O(n^2) ,
\]
where we omit the negligible contributions arising when there is just
one network geodesic between source and destination. We can sum these
contributions since the effect of double-counting is again negligible.
Thus, the total flow through $\origin$ is $\pi n^3 + O(n^2)$. If we
fix our attention on total flow through one of the bonds attached to
$\origin$
(so as to establish comparability with the results of Section \ref
{sec:first-order} for the Poissonian city), then we obtain $\frac{\pi
}{2} n^3 + O(n^2)$, compared with mean flow for the Poissonian city of $2n^3$.

However, account needs to be taken of the greater total length of the
grid network. Mean total network length produced by the unit intensity
Poisson line process is $\frac{\pi^2 n^2}{2}$, compared with
$2\times\pi n^2$ for the unit grid structure. Thus, a comparable grid
structure has to be based on segments of length $\frac{4}{\pi}$
rather than $1$. The flow produced through a bond attached to $\origin
$ by such a grid using the above protocol will be of order
\[
\frac{(\pi n^2)^2}{(\pi(({\pi}/{4})n)^2)^2} \times\frac{\pi
}{2} \biggl(\frac{\pi}{4}n\biggr)^3
=
2 n^3 ,
\]
where the second term computes the flow through a center bond when
rescaling from $n$ to $\frac{\pi}{4}n$, and the first factor is a
correction to ensure that total traffic is $(\pi n^2)^2$, not $(\pi
(\frac{\pi}{4}n)^2)^2$.
Thus, the Poissonian city is comparable to this grid structure in terms
of mean flow at the center. However, the grid structure with this
protocol can be shown to have the undesirable feature that,
asymptotically, a definite proportion of the total traffic (around 2\%)
occurs \textit{outside} the disc.

In fact, we can also carry out calculations for a somewhat more
demanding situation in which, asymptotically, the traffic stays within
the disc, as in the case of the Poissonian city; instead of supposing
the traffic to be divided equally between the two extreme network
geodesics, we suppose that it is divided equally among all possible
network geodesic connections. In effect, the actual network geodesic is
chosen uniformly at random, and in that case the probability that the
resulting network geodesic passes through the origin $(0,0)$ is
%
%
\begin{equation}\label{eq:grid-binomial}
\frac{
\mathbb{P}[\operatorname{Bin}(u+v,p)=u]
\mathbb{P}[\operatorname{Bin}(x+y,p)=x]
}{
\mathbb{P}[\operatorname{Bin}(u+v+x+y,p)=u+x]
}
=
\frac{{u+v\choose u}{x+y\choose x}}{{u+v+x+y\choose u+x}} .
\end{equation}
Here, $p=\frac12$, but the same result is in fact obtained for any
$0<p<1$. Choosing $p=\frac{u+x}{u+v+x+y}$, so as to control the
denominator, and applying Stirling's formula, we can use a Taylor
expansion to approximate this by
%
%
\begin{eqnarray}\label{eq:stirling}
&&\frac{1}{\sqrt{2\pi}}
\frac{(u+v+x+y)^{3/2}}{\sqrt{(u+v)(x+y)(u+x)(v+y)}}\nonumber\\[-8pt]\\[-8pt]
&&\qquad\times{}\exp\biggl(
-\frac{(u+v+x+y)(u y - x v)^2}{2(u+v)(x+y)(u+x)(v+y)}
\biggr) .\nonumber
\end{eqnarray}
Using polar coordinates based on an axis at $45^\circ$ to the Cartesian
axes and a Gaussian approximation based on $\sin^2(\theta-\phi
)\approx(\theta-\phi)^2$, we obtain a heuristic approximation for
large $n$ for the flow between two opposing quadrants of the disc:
%
%
\begin{equation}\label{eq:grid-asymptotic}
\mathop{\mathop{\sum\sum}_{u,v\geq0}}_{u^2+v^2\leq n^2}
\mathop{\mathop{\sum\sum}_{x,y\geq0}}_{x^2+y^2\leq n^2}
\frac{{u+v\choose u}{x+y\choose x}}{{u+v+x+y\choose u+x}}
\approx2n^3 .
\end{equation}
Conversion of this approach into a rigorous asymptotic argument would
require close attention to detailed asymptotics of the Binomial
distribution [\citet{Littlewood-1969}, \citet{McKay-1989}].
However, there
is an alternate argument which is more easily made rigorous: the
expression~(\ref{eq:grid-asymptotic}) can be reexpressed in terms of a
symmetric random walk $X$ as
\[
\mathop{\mathop{\sum\sum}_{u,v\geq0}}_{u^2+v^2\leq n^2}
\mathop{\mathop{\sum\sum}_{x,y\geq0}}_{x^2+y^2\leq n^2}
\mathbb{P}[X_{u+v}=v-u | X_{u+v+x+y}=v-u+y-x] .
\]
Under the $\mathbb{Z}$-action $u\to u-1, v\to v+1, x\to x+1, y\to
y+1$, a statistical pivot argument applied to the summand generates a
probability distribution on even integers or odd integers according to
the parity of $u+v$, $x+y$. An argument using the Hoeffding inequality
quantifies how this probability distribution concentrates around its
mode; the denominator is controlled by choosing the probability
$\mathbb{P}[X_{n+1}=X_n-1]=p=\frac{u+x}{u+v+x+y}$. It can thus be
shown that the
asymptotic behavior of the quadruple sum is given by the number of
$\mathbb{Z}$-orbits containing modal representatives close to the line
between $-(u,v)$ and $(x,y)$. This number can be expressed as a sum
susceptible to elementary asymptotic analysis, finally yielding a
rigorous argument for the asymptotic given in~(\ref{eq:grid-asymptotic}).

Consider now the total flow through a unit-length bond $\ell$
connected to the origin. This will equal half the total flow through
the origin, which itself can be viewed asymptotically as the sum of two
equal components from two different pairs of opposing quadrants. Thus,
this total flow is again asymptotic to $2n^3$.

Again, a grid structure comparable to the Poissonian city must be based
on segments of length $\frac{4}{\pi}$ rather than $1$, and a scaling
argument then shows that such a grid structure produces mean flow at
the center which is asymptotic to $\frac{4}{\pi}\times2
n^3=2.54648\ldots n^3$. Thus, traffic through the center under this
protocol is about $25\%$ higher in a comparable Manhattan city.

Geodesics in the Manhattan city are on average longer than those in the
Poisson line process; we can in fact argue in a manner analogous to
that of Section~\ref{sec:location-dependence} to show that mean
network distance between two independent uniformly random points in the
disc will be asymptotic to
\[
\frac{128}{45\pi} n \times \frac{1}{2\pi}\int_0^{2\pi}
(|{\sin\theta}|+|{\cos\theta}|) \,\dd\theta = \frac{128}{45\pi} n
\times\frac{4}{\pi} ,
\]
so the mean network flow over the whole disc for the grid will again be
about 25\% greater than mean network flow for the Poisson line process.

Of course, the second order behaviors of the flows are rather
different: flow at the center of the Poissonian city inherits
asymptotically nondegenerate randomness from the random configuration
of the Poisson line process, while a central limit argument shows that
the flows at the centers of the two kinds of flow in Manhattan cities
are asymptotically deterministic.

\section{Complements and conclusion}\label{sec:conclusion}

In conclusion, we present some notes concerning complements and issues
for further research, illustrating the potentially rich theory
concerning the Poissonian city.\vspace*{-1pt}

\subsection*{Empirical comparisons}
Clearly, the Poissonian city does not accurately represent real cities;
there will be variation both of geometry and of traffic flow. It would
be interesting to make empirical comparisons with actual street map and
traffic flow data, both in terms of network distance statistics
compared with the results of Section~\ref{sec:connections} and (much
more demanding from a data collection point of view) in terms of flow
statistics compared with the results of Section~\ref{sec:flow}. One
would expect qualitative agreement at best, rather than quantitative,
in view of the strong stochastic assumptions implicit in the Poissonian
city. Note, however, that the results on the asymptotic statistics of
flow at the center (Theorem~\ref{thm:scaling-limit}) reflect variation
across a sample of different cities, rather than within a particular city.\vspace*{-1pt}

\subsection*{Lower bounds on path length}
Compare the lower bound of Section~\ref{sec:true-geodesics}, Theorem
\ref{thm:poisson-lower-bound}, with the more general lower bound of
\citeauthor{AldousKendall-2007} [(\citeyear{AldousKendall-2007}),
Theorem~2], which holds for all connecting networks using total network
length proportional to $n$ based on patterns of nodes in a square
$[0,\sqrt{n}]^2$ satisfying a certain quantitative equidistribution
condition (related to an intuitive coupling construction). The
\citet{AldousKendall-2007} result provides an $\Omega(\sqrt{\log
n})$ bound on excess, whereas Theorem~\ref{thm:poisson-lower-bound}
uses detailed properties of Poisson line process networks to establish
a $\mathrm{const.}\times\log n$ lower bound. A~natural question is
whether there are any network constructions which provide
sub-logarithmic mean excess for appropriately equidistributed patterns
of nodes, or whether, on the other hand, the general lower bound can be
improved.\vspace*{-1pt}

%
%

\subsection*{Analytical characterization of limit}


The stochastic geometric construction of the limit distribution for
flow in the center of a Poissonian city using an improper stationary
anisotropic Poisson line process (Section~\ref{sec:limiting-object})
is explicit and lends itself to simulation; however, it would be
helpful also to have an analytic expression, or at least
characterization, of the limiting distribution. This seems difficult.
Note that we can produce a stochastic representation of the moment
generating function $M(p)$ of the limit distribution in terms of the
following probability: consider a Poisson process of intensity $\alpha
$ of \textit{pairs} of points on $[-1,1]\times(0,\infty)$. Then,
$M(-\alpha)$ is the probability that no pairs produced by this process
are separated from the origin by lines of the improper stationary
anisotropic Poisson line process. Similar representations are of use in
perfect simulation of area-interaction point process models [\citet
{Kendall-1997a}] and exact simulation of diffusions [\citet
{BeskosRoberts-2005a}]. However, in the current case it is not yet clear
whether this offers any progress toward simulation methods or
delivering a useful analytical representation.\vadjust{\goodbreak}

A slightly easier question is whether the convergence of $T_n/n^3$ to
the limit distribution holds for all moments. Again, at present no
progress on this can be offered beyond the work noted in Section \ref
{sec:uniform-integrability}.\vspace*{-1pt}

\subsection*{Aggregation issues}

What can we say about similar situations where the distribution of
nodes generating the flow is nonuniform? Or even when the nodes
generating the flows lie along the Poisson line process itself (thus
precluding the need for the ``cross country'' plumbing otherwise
required to get onto the network)? Considerations of this kind are
latent in the early work of \citet{Davidson-1974c}, and it would
be interesting to see them applied in the more quantitative setting of
the present work. It is possible that the coupling construction in
\citet{AldousKendall-2007} would be of use here.\vspace*{-1pt}

\subsection*{Three dimensions and higher}

In higher dimensions, one needs to consider what kind of network is
being deployed. One might, for example, consider the edge process of a
Poisson hyperplane tessellation or, alternatively, one might consider
network geodesics constrained to lie on the union of all faces of the
tessellation. In the second case, one can derive upper bounds on excess
by considering the derivative planar problem obtained by taking a
$2$-plane slice through the source and destination nodes; it is then a
question of how much the excess may be reduced by varying the
orientation of the slice, and it is a further question as to whether
the excess can be further substantially reduced by using paths which do
not lie wholly on the slicing plane. It may be possible to make
progress in the first case by adopting the growth process approach of
Section~\ref{sec:growth}.

Note that \citet{BoroczkySchneider-2008} describe
higher-dimensional results for similar problems; however, their results
concern standard stereological quantities, while we would need results
involving infima of lengths of regular curves on the boundaries of
Poisson cells.\vspace*{-1pt}

\subsection*{Moving beyond line processes}


Certainly, one can conceive of results for situations based on
processes which approximate Poisson line processes; Boolean models
based on long line segments or fiber processes for which there is
strong control of total fiber curvature. It would be particularly
interesting to determine the extent to which Poissonian cities and
Manhattan cities represent two extremes of a suitable class of models.\vspace*{-1pt}

\subsection*{User equilibrium}


The notion of UE [user equilibrium, \citet{Wardrop-1952},
contemporary with the related notion of Nash equilibrium] supposes that
each user has a utility structure for choosing which route they might
take based on travel time, which is affected not only by available
route lengths but also by flow along the routes. Interest is then
focused on systems of choices by users which result in user
equilibrium; no one user can obtain a shorter route by varying their
own route.
Explorations have already been made in the context of queueing
theory;\vadjust{\goodbreak}
see, for example, \citet{CalvertSolomonZiedins-1997}, who
consider the effect of augmenting a simple queueing network and
\citet{AfimeimoungaSolomonZiedins-2005}, who consider a system of
interactions between a $\cdot/M/1$ queue and a $\cdot/N^{(N)}/\infty
$ batch queue.
There are interesting possibilities in the context of the Poissonian
city, for example, considering that traffic from $p^-$ to $p^+$ chooses
each of the two possible routes prescribed by the semiperimeter
algorithm according to considerations both of length and of integrated
total flow along the routes.

Such problems are naturally formulated in terms of phase transitions in
statistical mechanics, perhaps using the improper stationary
anisotropic Poisson line process of Section~\ref{sec:limiting-object}.




\section*{Acknowledgments}
The author wishes to thank to Saul Jacka and Jon Warren for very
helpful conversations, and Ron Doney, Andreas Kyprianou, Juan Carlos
Pardo Millan and Mladen Savov for timely help with L\'evy process
theory.

\printaddresses


\begin{thebibliography}{40}

\bibitem[\protect\citeauthoryear{Afimeimounga, Solomon and
Ziedins}{2005}]{AfimeimoungaSolomonZiedins-2005}
%
\begin{barticle}[author]
\bauthor{\bsnm{Afimeimounga},~\bfnm{Heti}\binits{H.}},
\bauthor{\bsnm{Solomon},~\bfnm{Wiremu}\binits{W.}} \AND
\bauthor{\bsnm{Ziedins},~\bfnm{Ilze}\binits{I.}}
(\byear{2005}).
\btitle{The Downs--Thomson paradox: Existence, uniqueness and
stability of
user equilibria}.
\bjournal{Queueing Syst.}
\bvolume{49}
\bpages{321--334}.
\bmrnumber{MR2149647 (2006g:60134)}
\end{barticle}
%
\endbibitem

\bibitem[\protect\citeauthoryear{Aldous and
Kendall}{2008}]{AldousKendall-2007}
%
\begin{barticle}[author]
\bauthor{\bsnm{Aldous},~\bfnm{David~J.}\binits{D.~J.}} \AND
\bauthor{\bsnm{Kendall},~\bfnm{Wilfrid~S.}\binits{W.~S.}}
(\byear{2008}).
\btitle{Short-length routes in low-cost networks via Poisson line
patterns}.
\bjournal{Adv. in Appl. Probab.}
\bvolume{40}
\bpages{1--21}.
\bmrnumber{MR2411811 (2009f:60018)}
\end{barticle}
%
\endbibitem

\bibitem[\protect\citeauthoryear{Alsmeyer, Iksanov and
R{\oe}sler}{2009}]{AlsmeyerIksanovRosler-2009}
%
\begin{barticle}[author]
\bauthor{\bsnm{Alsmeyer},~\bfnm{Gerold}\binits{G.}},
\bauthor{\bsnm{Iksanov},~\bfnm{Alex}\binits{A.}} \AND
\bauthor{\bsnm{R{\oe}sler},~\bfnm{Uwe}\binits{U.}}
(\byear{2009}).
\btitle{On distributional properties of perpetuities}.
\bjournal{J. Theoret. Probab.}
\bvolume{22}
\bpages{666--682}.
\bmrnumber{MR2530108}
\end{barticle}
%
\endbibitem

\bibitem[\protect\citeauthoryear{Ambartzumian}{1990}]{Ambartzumian-1990}
%
\begin{bbook}[author]
\bauthor{\bsnm{Ambartzumian},~\bfnm{R.~V.}\binits{R.~V.}}
(\byear{1990}).
\btitle{Factorization Calculus and Geometric Probability}.
\bpublisher{Cambridge Univ. Press}, \baddress{Cambridge}.
\bmrnumber{MR1075011 (92b:60013)}
\end{bbook}
%
\endbibitem

\bibitem[\protect\citeauthoryear{Baccelli, Tchoumatchenko and
Zuyev}{2000}]{BaccelliTchoumatchenkoZuyev-2000}
%
\begin{barticle}[author]
\bauthor{\bsnm{Baccelli},~\bfnm{F.}\binits{F.}},
\bauthor{\bsnm{Tchoumatchenko},~\bfnm{K.}\binits{K.}} \AND
\bauthor{\bsnm{Zuyev},~\bfnm{S.}\binits{S.}}
(\byear{2000}).
\btitle{Markov paths on the Poisson--Delaunay graph with
applications to
routing in mobile networks}.
\bjournal{Adv. in Appl. Probab.}
\bvolume{32}
\bpages{1--18}.
\bmrnumber{MR1765174 (2001k:60010)}
\end{barticle}
%
\endbibitem

\bibitem[\protect\citeauthoryear{Baricz}{2008}]{Baricz-2008}
%
\begin{barticle}[author]
\bauthor{\bsnm{Baricz},~\bfnm{{\'A}rp{\'a}d}\binits{{\'A}.}}
(\byear{2008}).
\btitle{Mills' ratio: Monotonicity patterns and functional inequalities}.
\bjournal{J.~Math. Anal. Appl.}
\bvolume{340}
\bpages{1362--1370}.
\bmrnumber{MR2390935}
\end{barticle}
%
\endbibitem

\bibitem[\protect\citeauthoryear{Bertoin and Yor}{2001}]{BertoinYor-2001}
%
\begin{barticle}[author]
\bauthor{\bsnm{Bertoin},~\bfnm{Jean}\binits{J.}} \AND
\bauthor{\bsnm{Yor},~\bfnm{Marc}\binits{M.}}
(\byear{2001}).
\btitle{On subordinators, self-similar Markov processes and some
factorizations of the exponential variable}.
\bjournal{Electron. Comm. Probab.}
\bvolume{6}
\bpages{95--106 (electronic)}.
\bmrnumber{MR1871698 (2002k:60097)}
\end{barticle}
%
\endbibitem

\bibitem[\protect\citeauthoryear{Bertoin and Yor}{2002}]{BertoinYor-2002}
%
\begin{barticle}[author]
\bauthor{\bsnm{Bertoin},~\bfnm{Jean}\binits{J.}} \AND
\bauthor{\bsnm{Yor},~\bfnm{Marc}\binits{M.}}
(\byear{2002}).
\btitle{On the entire moments of self-similar Markov processes and
exponential functionals of L\'evy processes}.
\bjournal{Ann. Fac. Sci. Toulouse Math. (6)}
\bvolume{11}
\bpages{33--45}.
\bmrnumber{MR1986381 (2004d:60119)}
\end{barticle}
%
\endbibitem

\bibitem[\protect\citeauthoryear{Bertoin and Yor}{2005}]{BertoinYor-2005}
%
\begin{barticle}[author]
\bauthor{\bsnm{Bertoin},~\bfnm{Jean}\binits{J.}} \AND
\bauthor{\bsnm{Yor},~\bfnm{Marc}\binits{M.}}
(\byear{2005}).
\btitle{Exponential functionals of {L\'evy} processes}.
\bjournal{Probab. Surv.}
\bvolume{2}
\bpages{191--212 (electronic)}.
\bmrnumber{MR2178044 (2007b:60116)}
\end{barticle}
%
\endbibitem

\bibitem[\protect\citeauthoryear{Beskos and
Roberts}{2005}]{BeskosRoberts-2005a}
%
\begin{barticle}[author]
\bauthor{\bsnm{Beskos},~\bfnm{A.}\binits{A.}} \AND
\bauthor{\bsnm{Roberts},~\bfnm{G.~O.}\binits{G.~O.}}
(\byear{2005}).
\btitle{Exact simulation of diffusions}.
\bjournal{Ann. Appl. Probab.}
\bvolume{15}
\bpages{2422--2444}.
\bmrnumber{MR2187299 (2006e:60111)}
\end{barticle}
%
\endbibitem

\bibitem[\protect\citeauthoryear{Birnbaum}{1942}]{Birnbaum-1942}
%
\begin{barticle}[author]
\bauthor{\bsnm{Birnbaum},~\bfnm{Z.~W.}\binits{Z.~W.}}
(\byear{1942}).
\btitle{An inequality for Mill's ratio}.
\bjournal{Ann. Math. Statist.}
\bvolume{13}
\bpages{245--246}.
\bmrnumber{MR0006640 (4,19b)}
\end{barticle}
%
\endbibitem

\bibitem[\protect\citeauthoryear{B\"{o}r\"{o}czky and
Schneider}{2010}]{BoroczkySchneider-2008}
%
\begin{barticle}[author]
\bauthor{\bsnm{B\"{o}r\"{o}czky},~\bfnm{K\'{a}roly~J.}\binits
{K.~J.}} \AND
\bauthor{\bsnm{Schneider},~\bfnm{Rolf}\binits{R.}}
(\byear{2010}).
\btitle{The mean width of circumscribed random polytopes}.
\bjournal{Canad. Math. Bull.}
\bvolume{53}
\bpages{614--628}.
\end{barticle}
%
\endbibitem

\bibitem[\protect\citeauthoryear{Cabo and
Groeneboom}{1994}]{CaboGroeneboom-1994}
%
\begin{barticle}[author]
\bauthor{\bsnm{Cabo},~\bfnm{A.~J.}\binits{A.~J.}} \AND
\bauthor{\bsnm{Groeneboom},~\bfnm{P.}\binits{P.}}
(\byear{1994}).
\btitle{Limit theorems for functionals of convex hulls}.
\bjournal{Probab. Theory Related Fields}
\bvolume{100}
\bpages{31--55}.
\bmrnumber{MR1292189 (95g:60017)}
\end{barticle}
%
\endbibitem

\bibitem[\protect\citeauthoryear{Calvert, Solomon and
Ziedins}{1997}]{CalvertSolomonZiedins-1997}
%
\begin{barticle}[author]
\bauthor{\bsnm{Calvert},~\bfnm{Bruce}\binits{B.}},
\bauthor{\bsnm{Solomon},~\bfnm{Wiremu}\binits{W.}} \AND
\bauthor{\bsnm{Ziedins},~\bfnm{Ilze}\binits{I.}}
(\byear{1997}).
\btitle{Braess's paradox in a queueing network with state-dependent routing}.
\bjournal{J. Appl. Probab.}
\bvolume{34}
\bpages{134--154}.
\bmrnumber{MR1429062 (98f:60182)}
\end{barticle}
%
\endbibitem

\bibitem[\protect\citeauthoryear{Davidson}{1974}]{Davidson-1974c}
%
\begin{bincollection}[author]
\bauthor{\bsnm{Davidson},~\bfnm{Rollo}\binits{R.}}
(\byear{1974}).
\btitle{Line-processes, roads, and fibres}.
In \bbooktitle{Stochastic Geometry (A Tribute to the Memory of Rollo
Davidson)}
(\beditor{E. F. Harding and  D. G. Kendall}, eds.)
\bpages{248--251}.
\bpublisher{Wiley}, \baddress{London}.
\bmrnumber{MR0358975 (50 \#\#11431)}
\end{bincollection}
%
\endbibitem

\bibitem[\protect\citeauthoryear{Dufresne}{1990}]{Dufresne-1990}
%
\begin{barticle}[author]
\bauthor{\bsnm{Dufresne},~\bfnm{Daniel}\binits{D.}}
(\byear{1990}).
\btitle{The distribution of a perpetuity, with applications to risk theory
and pension funding}.
\bjournal{Scand. Actuar. J.}
\bvolume{1-2}
\bpages{39--79}.
\bmrnumber{MR1129194 (92i:62195)}
\end{barticle}
%
\endbibitem

\bibitem[\protect\citeauthoryear{Goldie and
Gr{\"u}bel}{1996}]{GoldieGrubel-1996}
%
\begin{barticle}[author]
\bauthor{\bsnm{Goldie},~\bfnm{Charles~M.}\binits{C.~M.}} \AND
\bauthor{\bsnm{Gr{\"u}bel},~\bfnm{Rudolf}\binits{R.}}
(\byear{1996}).
\btitle{Perpetuities with thin tails}.
\bjournal{Adv. in Appl. Probab.}
\bvolume{28}
\bpages{463--480}.
\bmrnumber{MR1387886 (97f:60124)}
\end{barticle}
%
\endbibitem

\bibitem[\protect\citeauthoryear{Groeneboom}{1988}]{Groeneboom-1988}
%
\begin{barticle}[author]
\bauthor{\bsnm{Groeneboom},~\bfnm{Piet}\binits{P.}}
(\byear{1988}).
\btitle{Limit theorems for convex hulls}.
\bjournal{Probab. Theory Related Fields}
\bvolume{79}
\bpages{327--368}.
\bmrnumber{MR959514 (89j:60024)}
\end{barticle}
%
\endbibitem

\bibitem[\protect\citeauthoryear{Hitczenko and
Weso{\l}owski}{2009}]{HitczenkoWesolowski-2010}
%
\begin{barticle}[author]
\bauthor{\bsnm{Hitczenko},~\bfnm{Pawe{\l}}\binits{P.}} \AND
\bauthor{\bsnm{Weso{\l}owski},~\bfnm{Jacek}\binits{J.}}
(\byear{2009}).
\btitle{Perpetuities with thin tails revisited}.
\bjournal{Ann. Appl. Probab.}
\bvolume{19}
\bpages{2080--2101}.
\bmrnumber{MR2588240}
\end{barticle}
%
\endbibitem

\bibitem[\protect\citeauthoryear{Kellerer}{1992}]{Kellerer-1992c}
%
\begin{bmisc}[author]
\bauthor{\bsnm{Kellerer},~\bfnm{H.~G.}\binits{H.~G.}}
(\byear{1992}).
\btitle{Ergodic behaviour of affine recursions III: Positive
recurrence and
null recurrence}.
\bnote{Technical report, Math. Inst. Univ. M\"unchen,
Theresienstrasse 39, 8000 M\"unchen, Germany}.
\end{bmisc}
%
\endbibitem

\bibitem[\protect\citeauthoryear{Kendall}{1997}]{Kendall-1997a}
%
\begin{bincollection}[vtex]
\bauthor{\bsnm{Kendall},~\bfnm{Wilfrid~S.}\binits{W.~S.}}
(\byear{1997}).
\btitle{On some weighted Boolean models}.
In \bbooktitle{Advances in Theory and Applications of Random Sets}
(\beditor{D. Jeulin}, ed.)
\bpages{105--120}.
\bpublisher{World Scientific}, \baddress{Singapore}.
\bmrnumber{MR1654418}
\end{bincollection}
%
\endbibitem

\bibitem[\protect\citeauthoryear{Kendall}{2008}]{Kendall-2009b}
%
\begin{barticle}[author]
\bauthor{\bsnm{Kendall},~\bfnm{Wilfrid~S.}\binits{W.~S.}}
(\byear{2008}).
\btitle{Networks and {Poisson} line patterns: Fluctuation asymptotics}.
\bjournal{Oberwolfach Rep.}
\bvolume{5}
\bpages{2670--2672}.
\end{barticle}
%
\endbibitem

\bibitem[\protect\citeauthoryear{Lamperti}{1972}]{Lamperti-1972}
%
\begin{barticle}[author]
\bauthor{\bsnm{Lamperti},~\bfnm{John}\binits{J.}}
(\byear{1972}).
\btitle{Semi-stable Markov processes. I}.
\bjournal{Z. Wahrsch. Verw. Gebiete}
\bvolume{22}
\bpages{205--225}.
\bmrnumber{MR0307358 (46 \#\#6478)}
\end{barticle}
%
\endbibitem

\bibitem[\protect\citeauthoryear{Littlewood}{1969}]{Littlewood-1969}
%
\begin{barticle}[author]
\bauthor{\bsnm{Littlewood},~\bfnm{J.~E.}\binits{J.~E.}}
(\byear{1969}).
\btitle{On the probability in the tail of a binomial distribution}.
\bjournal{Adv. in Appl. Probab.}
\bvolume{1}
\bpages{43--72}.
\bmrnumber{MR0240858 (39 \#\#2203)}
\end{barticle}
%
\endbibitem

\bibitem[\protect\citeauthoryear{McKay}{1989}]{McKay-1989}
%
\begin{barticle}[author]
\bauthor{\bsnm{McKay},~\bfnm{Brendan~D.}\binits{B.~D.}}
(\byear{1989}).
\btitle{On Littlewood's estimate for the binomial distribution}.
\bjournal{Adv. in Appl. Probab.}
\bvolume{21}
\bpages{475--478}.
\bmrnumber{MR997736 (90k:60030)}
\end{barticle}
%
\endbibitem

\bibitem[\protect\citeauthoryear{Miles}{1964}]{Miles-1964a}
%
\begin{barticle}[author]
\bauthor{\bsnm{Miles},~\bfnm{R.~E.}\binits{R.~E.}}
(\byear{1964}).
\btitle{Random polygons determined by random lines in a plane}.
\bjournal{Proc. Natl. Acad. Sci. USA}
\bvolume{52}
\bpages{901--907}.
\bmrnumber{MR0168000 (29 \#\#5265)}
\end{barticle}
%
\endbibitem

\bibitem[\protect\citeauthoryear{Narasimhan and
Smid}{2007}]{NarasimhanSmid-2007}
%
\begin{bbook}[author]
\bauthor{\bsnm{Narasimhan},~\bfnm{Giri}\binits{G.}} \AND
\bauthor{\bsnm{Smid},~\bfnm{Michiel}\binits{M.}}
(\byear{2007}).
\btitle{Geometric Spanner Networks}.
\bpublisher{Cambridge Univ. Press}, \baddress{Cambridge}.
\bmrnumber{MR2289615 (2009b:68002)}
\end{bbook}
%
\endbibitem

\bibitem[\protect\citeauthoryear{Pr{\"o}mel and
Steger}{2002}]{PromelSteger-2002}
%
\begin{bbook}[author]
\bauthor{\bsnm{Pr{\"o}mel},~\bfnm{Hans~J{\"u}rgen}\binits{H.~J.}}
\AND
\bauthor{\bsnm{Steger},~\bfnm{Angelika}\binits{A.}}
(\byear{2002}).
\btitle{The Steiner Tree Problem: A Tour Through Graphs, Algorithms, and Complexity}.
\bpublisher{Friedr. Vieweg \& Sohn}, \baddress{Braunschweig}.
\bmrnumber{MR1891564 (2003a:05047)}
\end{bbook}
%
\endbibitem

\bibitem[\protect\citeauthoryear{Rebolledo}{1980}]{Rebolledo-1980}
%
\begin{barticle}[author]
\bauthor{\bsnm{Rebolledo},~\bfnm{Rolando}\binits{R.}}
(\byear{1980}).
\btitle{Central limit theorems for local martingales}.
\bjournal{Z. Wahrsch. Verw. Gebiete}
\bvolume{51}
\bpages{269--286}.
\bmrnumber{MR566321 (81g:60023)}
\end{barticle}
%
\endbibitem

\bibitem[\protect\citeauthoryear{R{\'e}nyi and
Sulanke}{1968}]{RenyiSulanke-1968}
%
\begin{barticle}[author]
\bauthor{\bsnm{R{\'e}nyi},~\bfnm{A.}\binits{A.}} \AND
\bauthor{\bsnm{Sulanke},~\bfnm{R.}\binits{R.}}
(\byear{1968}).
\btitle{Zuf\"allige konvexe Polygone in einem Ringgebiet}.
\bjournal{Z. Wahrsch. Verw. Gebiete}
\bvolume{9}
\bpages{146--157}.
\bmrnumber{MR0229272 (37 \#\#4846)}
\end{barticle}
%
\endbibitem

\bibitem[\protect\citeauthoryear{Sampford}{1953}]{Sampford-1953}
%
\begin{barticle}[author]
\bauthor{\bsnm{Sampford},~\bfnm{M.~R.}\binits{M.~R.}}
(\byear{1953}).
\btitle{Some inequalities on Mill's ratio and related functions}.
\bjournal{Ann. Math. Statist.}
\bvolume{24}
\bpages{130--132}.
\bmrnumber{MR0054890 (14,995g)}
\end{barticle}
%
\endbibitem

\bibitem[\protect\citeauthoryear{Santal{\'o}}{1976}]{Santalo-1976}
%
\begin{bbook}[author]
\bauthor{\bsnm{Santal{\'o}},~\bfnm{Luis~A.}\binits{L.~A.}}
(\byear{1976}).
\btitle{Integral Geometry and Geometric Probability}.
\bpublisher{Addison-Wesley}, \baddress{Reading, MA}.
\bnote{With a foreword by Mark Kac, Encyclopedia of Mathematics and its
Applications, Vol. 1}.
\bmrnumber{MR0433364 (55 \#\#6340)}
\end{bbook}
%
\endbibitem

\bibitem[\protect\citeauthoryear{Steele}{1997}]{Steele-1997}
%
\begin{bbook}[author]
\bauthor{\bsnm{Steele},~\bfnm{J.~Michael}\binits{J.~M.}}
(\byear{1997}).
\btitle{Probability theory and combinatorial optimization}.
\bseries{CBMS-NSF Regional Conference Series in Applied Mathematics}
\bvolume{69}.
\bpublisher{SIAM},
\baddress{Philadelphia, PA}.
\bmrnumber{MR1422018 (99d:60002)}
\end{bbook}
%
\endbibitem

\bibitem[\protect\citeauthoryear{Stoyan, Kendall and
Mecke}{1995}]{StoyanKendallMecke-1995}
%
\begin{bbook}[author]
\bauthor{\bsnm{Stoyan},~\bfnm{D.}\binits{D.}},
\bauthor{\bsnm{Kendall},~\bfnm{Wilfrid~S.}\binits{W.~S.}} \AND
\bauthor{\bsnm{Mecke},~\bfnm{J.}\binits{J.}}
(\byear{1995}).
\btitle{Stochastic Geometry and Its Applications}, \bedition
{2nd} ed.
\bpublisher{Wiley}, \baddress{Chichester}.
\bnote{(First edition in 1987 joint with Akademie Verlag, Berlin.)}
\bmrnumber{MR895588 (88j:60034a)}
\end{bbook}
%
\endbibitem

\bibitem[\protect\citeauthoryear{Vervaat}{1979}]{Vervaat-1979}
%
\begin{barticle}[author]
\bauthor{\bsnm{Vervaat},~\bfnm{Wim}\binits{W.}}
(\byear{1979}).
\btitle{On a stochastic difference equation and a representation of
nonnegative infinitely divisible random variables}.
\bjournal{Adv. in Appl. Probab.}
\bvolume{11}
\bpages{750--783}.
\bmrnumber{MR544194 (81b:60064)}
\end{barticle}
%
\endbibitem

\bibitem[\protect\citeauthoryear{Voss, Gloaguen and
Schmidt}{2009}]{VossGloaguenSchmidt-2009}
%
\begin{bmisc}[author]
\bauthor{\bsnm{Voss},~\bfnm{F.}\binits{F.}},
\bauthor{\bsnm{Gloaguen},~\bfnm{C.}\binits{C.}} \AND
\bauthor{\bsnm{Schmidt},~\bfnm{V.}\binits{V.}}
(\byear{2009}).
\btitle{Scaling limits for shortest path lengths along the edges of
stationary tessellations}.
\bhowpublished{Preprint, Dept. Math, Univ. Ulm}.
\end{bmisc}
%
\endbibitem

\bibitem[\protect\citeauthoryear{Wardrop}{1952}]{Wardrop-1952}
%
\begin{barticle}[author]
\bauthor{\bsnm{Wardrop},~\bfnm{John~Glen}\binits{J.~G.}}
(\byear{1952}).
\btitle{Some theoretical aspects of road traffic research}.
\bjournal{Proceedings, Institute of Civil Engineers, Part II}
\bvolume{1}
\bpages{325--378}.
\end{barticle}
%
\endbibitem

\bibitem[\protect\citeauthoryear{Whitt}{2007}]{Whitt-2007}
%
\begin{barticle}[author]
\bauthor{\bsnm{Whitt},~\bfnm{Ward}\binits{W.}}
(\byear{2007}).
\btitle{Proofs of the martingale FCLT}.
\bjournal{Probab. Surv.}
\bvolume{4}
\bpages{268--302}.
\bmrnumber{MR2368952 (2008k:60079)}
\end{barticle}
%
\endbibitem

\bibitem[\protect\citeauthoryear{Yor}{1992}]{Yor-1992a}
%
\begin{barticle}[author]
\bauthor{\bsnm{Yor},~\bfnm{Marc}\binits{M.}}
(\byear{1992}).
\btitle{On some exponential functionals of Brownian motion}.
\bjournal{Adv. in Appl. Probab.}
\bvolume{24}
\bpages{509--531}.
\bmrnumber{MR1174378 (94b:60095)}
\end{barticle}
%
\endbibitem

\bibitem[\protect\citeauthoryear{Yukich}{1998}]{Yukich-1998}
%
\begin{bbook}[author]
\bauthor{\bsnm{Yukich},~\bfnm{Joseph~E.}\binits{J.~E.}}
(\byear{1998}).
\btitle{Probability Theory of Classical Euclidean Optimization Problems}.
\bseries{Lecture Notes in Math.}
\bvolume{1675}.
\bpublisher{Springer}, \baddress{Berlin}.
\bmrnumber{MR1632875 (2000d:60018)}
\end{bbook}
%
\endbibitem

\end{thebibliography}
\end{document}